\documentclass[nospthms,envcountsect]{svmult} 
\usepackage{makeidx}
\makeindex
\usepackage{latexsym}
\usepackage{amsfonts}
\usepackage{amsmath}
\usepackage{amssymb}
\usepackage{amscd}
\usepackage{multicol}    
\usepackage{enumerate}
    \newtheorem{theorem}                    {Theorem}       [section]
    \newtheorem{lemma}      [theorem]       {Lemma}
    \newtheorem{corollary}  [theorem]       {Corollary}
    \newtheorem{proposition}[theorem]       {Proposition}

\begin{document}
\catcode`@=11
\atdef@ I#1I#2I{\CD@check{I..I..I}{\llap{$\m@th\vcenter{\hbox
  {$\scriptstyle#1$}}$}
  \rlap{$\m@th\vcenter{\hbox{$\scriptstyle#2$}}$}&&}}
\atdef@ E#1E#2E{\ampersand@
  \ifCD@ \global\bigaw@\minCDarrowwidth \else \global\bigaw@\minaw@ \fi
  \setboxz@h{$\m@th\scriptstyle\;\;{#1}\;$}%
  \ifdim\wdz@>\bigaw@ \global\bigaw@\wdz@ \fi
  \@ifnotempty{#2}{\setbox@ne\hbox{$\m@th\scriptstyle\;\;{#2}\;$}%
    \ifdim\wd@ne>\bigaw@ \global\bigaw@\wd@ne \fi}%
  \ifCD@\enskip\fi
    \mathrel{\mathop{\hbox to\bigaw@{}}%
      \limits^{#1}\@ifnotempty{#2}{_{#2}}}%
  \ifCD@\enskip\fi \ampersand@}
\catcode`@=\active

\renewcommand{\labelenumi}{\alph{enumi})}
\newcommand{\chr}{\operatorname{char}}
\newcommand{\isom}{\stackrel{\sim}{\longrightarrow}}
\newcommand{\Aut}{\operatorname{Aut}}
\newcommand{\Hom}{\operatorname{Hom}}
\newcommand{\End}{\operatorname{End}}
\newcommand{\HOM}{\operatorname{{\mathcal H{\frak{om}}}}}
\newcommand{\EXT}{\operatorname{\mathcal E{\frak xt}}}
\newcommand{\Tot}{\operatorname{Tot}}
\newcommand{\Ext}{\operatorname{Ext}}
\newcommand{\Gal}{\operatorname{Gal}}
\newcommand{\Pic}{\operatorname{Pic}}
\newcommand{\Spec}{\operatorname{Spec}}
\newcommand{\trdeg}{\operatorname{trdeg}}
\newcommand{\im}{\operatorname{im}}
\newcommand{\coim}{\operatorname{coim}}
\newcommand{\coker}{\operatorname{coker}}
\newcommand{\gr}{\operatorname{gr}}
\newcommand{\id}{\operatorname{id}}
\newcommand{\Br}{\operatorname{Br}}
\newcommand{\cd}{\operatorname{cd}}
\newcommand{\CH}{CH}
\newcommand{\Alb}{\operatorname{Alb}}
\renewcommand{\lim}{\operatornamewithlimits{lim}}
\newcommand{\colim}{\operatornamewithlimits{colim}}
\newcommand{\rk}{\operatorname{rank}}
\newcommand{\codim}{\operatorname{codim}}
\newcommand{\NS}{\operatorname{NS}}
\newcommand{\cone}{{\rm cone}}
\newcommand{\rank}{\operatorname{rank}}
\newcommand{\ord}{{\rm ord}}
\newcommand{\f}{{\cal G}}
\newcommand{\g}{{\cal F}}
\newcommand{\du}{{\cal D}}
\newcommand{\G}{{\mathbb G}}
\newcommand{\N}{{\mathbb N}}
\newcommand{\A}{{\mathbb A}}
\newcommand{\Z}{{{\mathbb Z}}}
\newcommand{\Q}{{{\mathbb Q}}}
\newcommand{\R}{{{\mathbb R}}}
\newcommand{\B}{{\mathbb Z}^c}
\renewcommand{\H}{{{\mathbb H}}}
\renewcommand{\P}{{{\mathbb P}}}
\newcommand{\F}{{{\mathbb F}}}
\newcommand{\m}{{\mathfrak m}}
\newcommand{\Sch}{{\text{\rm Sch}}}
\newcommand{\et}{{\text{\rm et}}}
\newcommand{\Zar}{{\text{\rm Zar}}}
\newcommand{\Nis}{{\text{\rm Nis}}}
\newcommand{\tr}{\operatorname{tr}}
\newcommand{\tor}{{\text{\rm tor}}}
\newcommand{\red}{{\text{\rm red}}}
\newcommand{\Div}{\operatorname{Div}}
\newcommand{\Ab}{{\text{\rm Ab}}}
\newcommand{\DD}{{\mathbb Z}^c}
\renewcommand{\div}{\operatorname{div}}
\newcommand{\corank}{\operatorname{corank}}
\renewcommand{\O}{{\cal O}}
\newcommand{\C}{{\mathbb C}}
\newcommand{\p}{{\mathfrak p}}
\newcommand{\proof}{\noindent{\it Proof. }}
\newcommand{\proofend}{\hfill $\Box$ \\}
\newcommand{\rem}{\noindent {\it Remark. }}
\newcommand{\example}{\noindent {\bf Example. }}
\newcommand{\ar}{{\text{\rm ar}}}
\newcommand{\del}{{\delta}}
\title*{Duality via cycle complexes}
\author{Thomas Geisser\thanks{Supported in part by NSF grant No.0556263}}
\institute{University of Southern California}

\maketitle

\vskip-1.5cm

\hfill{Dedicated to the memory of}

\hfill{Hermann Braun}

\begin{abstract}
We show that Bloch's complex of relative zero-cycles can
be used as a dualizing complex over perfect fields and number rings.
This leads to duality theorems for torsion sheaves on arbitrary
separated schemes of finite type over algebraically closed fields,
finite fields, local fields of mixed characteristic, and rings
of integers in number rings, generalizing results which so far
have only been known for smooth schemes or in low dimensions,
and unifying the $p$-adic and $l$-adic theory.
As an application, we generalize Rojtman's theorem to normal,
projective schemes.
\end{abstract}

\section{Introduction}
If $f:X\to S$ is separated and of finite type, then in order to
obtain duality theorems from the adjointness
$$ R\Hom_X(\f,Rf^!\g)\cong R\Hom_S(Rf_!\f,\g)$$
for torsion \'etale sheaves $\f$ on $X$ and $\g$ on $S$,
one has to identify the complex $Rf^!\g$. For example, if $f$
is smooth of relative dimension $d$ and if $m$ is invertible
on $S$, then Poincar\'e duality of SGA 4 XVIII states that
$Rf^!\g\cong f^*\g\otimes\mu_m^{\otimes d}[2d]$
for $m$-torsion sheaves $\g$. We show that if $S$ is the spectrum
of a perfect field,
or a Dedekind ring of characteristic $0$ with perfect residue fields, then
Bloch's complex of zero-cycles can by used to explicitly
calculate $Rf^!\g$.

For a scheme $X$ essentially of finite type over $S$,
let $\DD_X$ be the complex of \'etale sheaves which in degree $-i$
associates to $U\to X$ the free abelian group generated by cycles
of relative dimension $i$ over $S$ on $U\times_S\Delta^i$
which meet all faces properly, and alternating sum of intersection
with faces as differentials \cite{bloch}. Over a field, the
higher Chow group of zero cycles $CH_0(X,i)$ is by definition the $-i$th
cohomology of the global sections $\DD_X(X)$, and if
$X$ is smooth of relative dimension $d$ over a perfect field $k$
of characteristic $p$, then $\DD_X/m\cong \mu_m^{\otimes d}[2d]$ for
$m$ prime to $p$, and $\DD_X/p^r\cong \nu_r^d[d]$,
the logarithmic de Rham-Witt sheaf \cite{bkbl, marcI}.
Our main result is that if $f:X\to Y$ is
separated and of finite type over a perfect field $k$,
and if $\g$ is a torsion sheaf on $X$, then there is a quasi-isomorphism
\begin{equation}\label{koko36}
R\Hom_X(\g,\DD_X)\cong R\Hom_Y(Rf_!\g,\DD_Y).
\end{equation}
If $k$ is algebraically closed and
$\g$ is constructible, this yields perfect pairings of finite groups
$$\Ext^{1-i}_X(\g,\DD_X)\times H^i_c(X_\et,\g)\to \Q/\Z.$$
In particular, we obtain an isomorphism
$CH_0(X,i,\Z/m)\cong H^i_c(X_\et,\Z/m)^*$ of finite groups,
generalizing a theorem of Suslin \cite{suslinetale} to arbitrary $m$,
and an isomorphism of the abelianized (profinite) fundamental group
$\pi_1^{ab}(X)$ with $CH_0(X,1,\hat \Z)$ for any proper scheme $X$
over $k$.
As an application, we generalize Rojtman's theorem
\cite{milneroitman, roitman} to normal schemes $X$,
projective over an algebraically closed field: The Albanese map
induces an isomorphism
$${}_\tor CH_0(X)\cong {}_\tor \Alb_X(k)  $$
between the torsion points of the Chow group of zero-cycles on $X$
and the torsion points of the Albanese variety in the sense of Serre
\cite{serre}. This is a homological version of
Rojtman's theorem which differs from the cohomological version
of Levine and Krishna-Srinivas \cite{levinerojtman, srin} relating the
Albanese variety in the sense of Lang-Weil to the Chow group defined by
Levine-Weibel \cite{levineweibel}. We give an example to show that
for non-normal schemes, the torsion elements of $CH_0(X)$
cannot be parametrized by an abelian variety in general.

If $k$ is finite, and $X$ is
separated and of finite type over $k$, we obtain
for constructible sheaves $\g$ perfect pairings of finite groups
$$\Ext^{2-i}_X(\g,\DD_X)\times H^i_c(X_\et,\g)\to \Q/\Z.$$
This generalizes results of Deninger \cite{deninger2} for curves,
Spie\ss\ \cite{spiess} for surfaces, and Milne \cite{milnevalues} and
Moser \cite{moser} for the $p$-part in characteristic $p$.
In fact, the dualizing complex of Deninger is
quasi-isomorphic to $\DD_X$ by a result of Nart \cite{nart},
and the  niveau spectral sequence of $\DD_X/p^r$ degenerates to
the dualizing complex of Moser. If $\pi_1^{ab}(X)^0$ is the kernel of
$\pi_1^{ab}(X)\to \Gal(k)$,
then for $X$ proper and $\bar X=X\times_k\bar k$,
we obtain a short exact sequence
$$0\to CH_0(\bar X,1)^\wedge_{G}\to  \pi_1^{ab}(X)^0
\to CH_0(\bar X)^G\to 0.$$
We obtain a similar duality theorem
for schemes over local fields of characteristic $0$.

If $f:X\to S$ is a scheme over the spectrum of a Dedekind ring of
characteristic $0$ with perfect residue fields, then, assuming the
Beilinson-Lichtenbaum
conjecture, there is a quasi-isomorphism
\begin{equation}\label{thom}
R\Hom_X(\g,\DD_X)\cong R\Hom_S(Rf_!\g,\DD_S).
\end{equation}
Even though the duality theorem over a field of characteristic $p$
could in principle be
formulated by treating the prime to $p$-part and $p$-part separately,
it is clear that over a Dedekind ring one needs a complex which
treats both cases uniformly.

If $S$ is the ring of integers in a number ring,
and if we define cohomology with compact support
$H^i_c(X_\et,\g)$ as the cohomology of the complex $R\Gamma_c(S_\et,Rf_!\g)$,
where $R\Gamma_c(S_\et,-)$ is cohomology with compact support of $S$
\cite{mazur, adt}, then combining \eqref{thom} with Artin-Verdier duality,
we get perfect pairings of finite groups
$$ \Ext^{2-i}_X(\g,\DD_X)\times H^i_c(X_\et,\g)\to \Q/\Z$$
for constructible $\g$.
This generalizes results of Artin-Verdier \cite{mazur}
for $\dim X=1$, Milne \cite{adt} for $\dim X=1$ and
$X$ possibly singular, or $X$ smooth over $S$, and Spie\ss\
\cite{spiess} for $\dim X=2$.

If $S$ is a henselian discrete valuation
ring of mixed characteristic with closed point $i:s\to S$,
and if we define cohomology with compact support in the closed
fiber $H^i_{X_s,c}(X_\et,\g)$ to be the cohomology of
$R\Gamma(S_\et,i_*Ri^!Rf_!\g)$, then there are perfect pairings
of finite groups
$$ \Ext^{2-i}_X(\g,\DD_X)\times H^i_{X_s,c}(X_\et,\g)\to \Q/\Z$$
for constructible $\g$.

We outline the proof of our main theorem \eqref{koko36}.
The key observation is that for $i:Z\to X$ a closed embedding
over a perfect field, we have a quasi-isomorphism of complexes of \'etale
sheaves $Ri^!\DD_X\cong \DD_Z$ (purity).
In order to prove purity, we show that $\DD_X$ has
\'etale hypercohomological descent over algebraically closed fields
(i.e. its cohomology and \'etale hypercohomology agree), and then use purity
for the cohomology of $\DD_X$ proved by Bloch \cite{blocloc} and Levine
\cite{levinemoving}. To prove \'etale hypercohomological descent, 
we use the argument
of Thomason \cite{thomason} to reduce to finitely generated fields over $k$,
and in this case use results of Suslin \cite{suslinetale} for the prime
to $p$-part, and Geisser-Levine \cite{marcI} and
Bloch-Kato \cite{blochkato} for the $p$-part.

Having purity, an induction and devissage argument is used to reduce to
the case of a constant sheaf on a smooth and proper scheme,
in which case we check that our pairing agrees with the classical
pairing of SGA 4 XVIII and Milne \cite{milnevalues} for the
prime to $p$ and $p$-primary part, respectively.

Throughout the paper, scheme over $S$ denotes a
separated scheme of finite type over $S$. We always work on the small
\'etale site of a scheme. For an abelian group $A$,
we denote by $A^*=\Hom(A,\Q/\Z)$ its Pontrjagin dual, by
$A^\wedge =\lim A/m$ its pro-finite completion, by ${}_mA$ the $m$-torsion
of $A$, and by $TA=\lim {}_mA$ its Tate-module.

{\it Acknowledgements.}
The work of this paper was inspired by the work of, and discussions
with, U.Jannsen, S.Lichtenbaum, S.Saito and K.Sato.
We are indebted to the referee for his careful reading and
helpful suggestions.

\section{The dualizing complex}
We recall some properties of Bloch's higher Chow complex \cite{bloch}, see
\cite{handbook} for a survey and references.
For a fixed regular scheme $S$ of finite Krull dimension $d$,
and an integral scheme $X$ essentially
of finite type over $S$, the relative dimension of $X$ over $S$ is
$$ \dim_S X=\trdeg(k(X):k(\p)) -\text{ht}\; \p +d.$$
Here $k(X)$ is the function field of $X$,
$\p$ is the image of the generic point of $X$ in $S$,
and $k(\p)$ its residue field.
For example, if $K$ is an extension field of transcendence degree $r$ over $k$,
then the dimension of $K$ over $k$ is $r$.
For a scheme $X$ essentially of finite type over $S$, we
define $z_n(X,i)$ to be the free abelian group generated by closed
integral subschemes of relative dimension $n+i$ over $S$ on
$X\times_S\Delta^i$ which meet all
faces properly. If $z_n(X,*)$ is the complex of abelian groups
obtained by taking the alternating sum of intersection with face
maps as differentials, then $z_n(-,*)$ is a (homological) complex
of sheaves for the \'etale topology \cite{ichdede}.
We define $\DD_X(n)=z_n(-,*)[2n]$ to be the (cohomological) complex with
the \'etale sheaf $z_n(-,-i-2n)$ in degree $i$.
If $n=0$, then we sometimes write $\DD_X$ instead of $\DD_X(0)$,
and we sometimes omit $X$ if there is no ambiguity.
For a quasi-finite, flat map $f:X\to Y$, we have a pull-back
$f^*\DD_Y(n)\to \DD_X(n)$ because $z_n(-,*)$ is contravariant for
such maps, and for a proper map $f:X\to Y$ we have a push-forward
$f_*\DD_X(n)\to \DD_Y(n)$ because $z_n(-,*)$ is covariant for
such maps. We frequently use that $\DD_X(n)$ is a complex of flat
sheaves, hence tensor product and derived tensor product with
$\DD_X(n)$ agree.

From now on we assume that the base $S$ is the
spectrum of a field or of a Dedekind ring. Then
for a closed embedding $i:Z\to X$ over $S$, we have a quasi-isomorphism
$\DD_Z(n)\cong Ri^!\DD_X(n)$
on the Zariski-site \cite{blocloc}, \cite{levinemoving}. We will refer
to this fact as purity or the localization property; the fact that the
analog statement holds on the \'etale site if $n\leq 0$ and the residue
fields of $S$ are perfect is a key result of this paper.
If $p:X\times {\Bbb A}^r\to X$ is the projection, we have a
quasi-isomorphism of complexes of Zariski-sheaves
\begin{equation}\label{affineform}
p_*\DD_{X\times {\Bbb A}^r}(n)\cong \DD_X(n-r)[2r].
\end{equation}
For a Grothendieck topology $t$, we define
$$H_i(X_t,\DD(n))=H^{-i}R\Gamma(X_t,\DD(n)).$$
We sometimes omit the $t$ when we use the Zariski topology, and note that
for $X$ over a field or discrete valuation ring \cite{bloch},
\cite{levinemoving}
$$CH_n(X,i-2n)cong H_i(X_\Zar,\DD(n)).$$
If $\Z_X(n)$ is the motivic complex of Voevodsky \cite{voevodsky}, then
on a smooth scheme $X$ of dimension $d$ over a field $k$,
\begin{equation}\label{pdual}
\DD_X(n)\cong \Z_X(d-n)[2d].
\end{equation}
The Beilinson-Lichtenbaum conjecture (translated into homological
notation by \eqref{pdual}) states that for a scheme $X$ as above, and
$m$ prime to the characteristic of $k$, the change of topology map
induces an isomorphism
$$H_i(X_\Zar,\DD_X/m(n))\xrightarrow{\sim} H_i(X_\et,\DD_X/m(n))$$
for $i\geq d+n$.
If $m$ is a power of the characteristic of $p$, then the analog
statement is known \cite{marcI}, and with $\Q$-coeffients,
Zariski and \'etale hypercohomology agree for smooth schemes and
all $i$. The Beilinson-Lichtenbaum conjecture is implied by the
Bloch-Kato conjecture stating that for any field $F$ of characteristic
prime to $m$, the norm residue homomorphism
$K_n^M(F)/m\to H^n(F_\et,\mu_m^{\otimes n})$ is an isomorphism for all $n$
\cite{bkbl}. A proof of the Bloch-Kato conjecture is announced by
Rost and Voevodsky, but since there is no published account at this
time, we will point out any use of the Beilinson-Lichtenbaum conjecture.

We define $H_i(F_t,\DD(n))=\colim_U H_i(U_t,\DD(n))$ for $F$ a finitely
generated field $F$ over $S$, where the
colimit runs through $U$ of finite type over $S$ with field of
functions $F$. Then $H_i(F_t,\DD(n))\cong H^{2d-i}(F_t,\Z(d-n))$ if
$F$ has dimension $d$ over $S$ by \eqref{pdual}. In particular, this group
vanishes for the Zariski topology if $i<d+n$, and agrees with
$K_{d-n}^M(F)$ for $i=d+n$. The Beilinson-Lichtenbaum conjecture
implies that $H_i(F,\Z^c(n))\cong H_i(F_\et,\DD(n))$ for $i\geq d+n$.

\begin{proposition}\label{zarnivlemma}
Let $X$ be a scheme over a field or Dedekind ring. Then
there are spectral sequences
\begin{equation}\label{niveau}
E^1_{s,t}=\bigoplus_{x\in X_{(s)}}H^{s-t}(k(x),\Z(s-n)) \Rightarrow
H_{s+t}(X,\DD(n)).
\end{equation}
In particular, $H_i(X,\DD(n))=0$ for $i<n$.
\end{proposition}

\proof
If we let $F_s \DD(n)$ be the subcomplex generated by cycles
of dimension $n+i$ on $X\times \Delta^i$ such that the projection to
$X$ has dimension at most $s$ over $S$, then we get the spectral sequence
$$E^1_{s,t}=H_{s+t}(X,F_s/F_{s-1}\DD(n))
\Rightarrow H_{s+t}(X,\DD(n)).$$
By the localization property, we get
\begin{multline*}
$$H_{s+t}(X,F_s/F_{s-1}\DD(n))
\cong \displaystyle\bigoplus_{x\in X_{(s)}}H_{s+t}(k(x),\DD(n))\\
\cong \bigoplus_{x\in X_{(s)}}H^{s-t}(k(x),\Z(s-n)).
\end{multline*}
An inspection shows that $E_{s,t}^1=0$ for $t<n$, hence the
vanishing. \proofend

Recall that $\nu_r^i=W_r\Omega^i_{X,log}$ is the logarithmic de Rham-Witt
sheaf. The following Proposition is a a finite coefficient-version
of \eqref{pdual}:

\begin{proposition}\label{identifyme}
Let $X$ is smooth of dimension $d$ over a perfect field $k$, and 
let $n\leq d$. Then there are quasi-isomorphisms of complexes of 
\'etale sheaves
\begin{align*}
\DD_X/m(n)&\cong \mu_m^{\otimes d-n}[2d] \qquad\;\;\text{if } \chr k\not|m; \\
\DD_X/p^r(n)& \cong \nu^{d-n}_r[d+n] \qquad \text{if $p=\chr k$}.
\end{align*}
This is compatible with the Gysin maps
$H_j(Z_\et,\DD_Z(n))\to H_j(X_\et,\DD_X(n))$, for closed embeddings
$i:Z\to X$ of pure codimension between smooth schemes.
\end{proposition}

\proof
The prime to $p$-part has been proved in
\cite[Thm. 4.14, Prop.4.5(2)]{bkbl}. For the $p$-primary part,
the quasi-isomorphism is given by \cite{marcI}
\begin{multline*}
\DD/p^r(n)[-d-n]\xleftarrow{\sim}\tau_{\leq 0}(\DD/p^r(n)[-d-n])
\xrightarrow{\sim}{\cal H}^{-d-n}(\DD/p^r(n))\\
\xrightarrow{\sim} G({\cal H}^{-d-n}(\DD/p^r(n)))
\cong G(\nu^{d-n}_r)\xleftarrow{\sim} \nu^{d-n}_r,
\end{multline*}
where $G(\nu^{d-n}_r)$ is the Gersten resolution arising as the
$E^1$-complex of the niveau spectral sequence \eqref{niveau}
$$ \bigoplus_{x\in X^{(0)}}(i_x)_*\nu_{r,k(x)}^{d-n}\to
\bigoplus_{x\in X^{(1)}}(i_x)_*\nu_{r,k(x)}^{d-n-1}\to\cdots, $$
and similarly for ${\cal H}^{-d-n}(\DD/p^r(n))$.
The Gersten resolutions identify via the isomorphisms
$H_{d+n-i}(k(x),\Z^c/p^r(n)) \xleftarrow{\sim}
K_{d-n-i}^M(k(x))/p^r\xrightarrow{\sim}
\nu_{r,k(x)}^{d-n-i}$ for a field $k(x)$ of codimension $i$,
i.e. transcendence degree $d-i$.
The compatibility of cohomology with proper push-forward and
flat equidimensional pull-back follows from the corresponding
property of the Gersten resolution \cite[Prop.1.18]{grossuwa}.
\proofend

The $p$-primary part of Proposition \ref{identifyme} can
be generalized to singular schemes.
Let $\tilde\nu_r(0)$ be the complex of \'etale sheaves
$$ \cdots \to \bigoplus_{x\in X_{(1)}}(i_x)_*\nu_{r,k(x)}^1\to
\bigoplus_{x\in X_{(0)}}(i_x)_*\nu_{r,k(x)}^0\to 0$$ used by Moser
\cite[1.5]{moser} (loc.cit. indexes by codimension, which makes the
treatment more complicated).

\begin{proposition}\label{locokp}
Let $X$ be a separated scheme of finite type over a perfect field $k$
of characteristic $p$. Then ${\B}/p(n)\cong 0$ for $n<0$, and there are 
isomorphisms of \'etale sheaves 
${\mathcal H}_i({\B}/p^r(0))\cong {\mathcal H}_i(\tilde\nu_r(0))$, 
compatible with proper push-forward.
In particular,
$H_i(X_\et,{\B}/p^r(0))\cong H_i(X_\et,\tilde\nu_r(0))$.
\end{proposition}

\proof
If $n<0$, let $R$ be a finitely generated algebra over $k$.
Write $R$ as a quotient of a smooth algebra $A$, and let $U=\Spec
A-\Spec R$. Then the localization sequence for higher Chow groups
$$ \cdots \to H_{i+1}(U,{\B}/p(n))\to H_i(\Spec R,{\B}/p(n))\to
H_i(\Spec A,{\B}/p(n)) \to \cdots $$ together with the fact that
$H_i(X_\Zar,{\B}/p(n))=0$ for smooth $X$ and $n<0$ \cite{marcI}
shows that $H_i(\Spec R,{\B}/p(n))\cong H_i({\B}/p(n)(\Spec R))\cong 0$. For
$n=0$, we can assume that $k$ is algebraically closed.
Consider the spectral sequence \eqref{niveau} with mod $p^r$-coefficients.
The $E^1_{s,t}$-terms vanish for $t<0$, and according to \cite{marcI}, they
also vanish for $t>0$.
Since $H^s(k(x),\Z/p^r(s))\cong H^0(k(x)_\et,\nu_r^s)$,
the cohomology of $\DD/p^r(0)(X)$ agrees with the cohomology of the
complex $\tilde\nu_r(0)(X)$ in a functorial way.
\proofend

It would be interesting to write down a map
$\DD/p^r(0)\to \tilde\nu_r(0)$ of complexes inducing
the isomorphism on cohomology of Proposition \ref{locokp}.
We will see
below that there is a quasi-isomorphism
$\DD/p^r(0)\xrightarrow{\sim} Rf^!\Z/p^r $. On the other hand,
Jannsen-Saito-Sato \cite{JSS} show that there is a quasi-isomorphism
$ \tilde \nu_r(0)\xrightarrow{\sim} Rf^!\Z/p^r $.

\begin{lemma}\label{intmodm}
If $\g$ is an $m$-torsion sheaf, then we have a quasi-isomorphism
$$ R\Hom_{X,\Z/m}(\g,\DD_X/m)[-1]\cong R\Hom_X(\g,\DD_X).$$
\end{lemma}

\proof
The exact, fully faithful inclusion functor
$F:Shv_{\Z/m}(X)\to Shv_{\Z}(X)$ from \'etale sheaves of $\Z/m$-modules
to \'etale sheaves of abelian groups has the left adjoint
$-\otimes_\Z\Z/m$ and the right adjoint "$m$-torsion" $T_m=\Hom_\Z(\Z/m,-)$;
in particular, $T_m$ preserves injectives. Moreover,
the left derived functor $-\otimes_\Z^L\Z/m$
of the tensor product agrees with the shift $RT_m[1]$
of the right derived functor of the $m$-torsion functor as functors
$D(Shv_{\Z}(X))\to D(Shv_{\Z/m}(X))$.
Indeed, both are quasi-isomorphic to the double complex
$C^\cdot\xrightarrow{\cdot m}C^\cdot$.
Since $\DD_X$ consists of flat sheaves, we have
$\DD_X/m\cong \DD_X\otimes^L\Z/m \cong RT_m\DD_X[1]$.
The Lemma follows from
$$ R\Hom_X(F\g,\DD_X)\cong R\Hom_{X,\Z/m}(\g,RT_m\DD_X)
\cong R\Hom_{X,\Z/m}(\g,\DD_X/m[-1])$$
\proofend

We also use frequently that for a complex of torsion abelian groups $C^\cdot$
we have
$\Hom(C^\cdot,\Q/\Z)[-1]\cong R\Hom(C^\cdot,\Z)$, and in particular
$H^iR\Hom(C^\cdot,\Z)\cong H^{1-i}(C^\cdot)^*$.

\section{Etale descent}
The main result of this section is purity and a trace map for
$\DD(0)$. We give a conceptual proof assuming the Beilinson-Lichtenbaum
conjecture, and an ad-hoc proof of a weaker, but for our purposes
sufficient result, avoiding the use of the Beilinson-Lichtenbaum conjecture.
We first use an argument of Thomason \cite{thomason} to show that for
$n\leq 0$, $\DD(n)$ has \'etale cohomological descent.

\begin{theorem}\label{locok}
Assume the Beilinson-Lichtenbaum conjecture holds for schemes
over the algebraically closed field $k$.
If $X$ is a scheme over $k$ and $n\leq 0$, then
$$\DD(n)(X)\cong R\Gamma(X_\et,\DD(n)).$$
\end{theorem}

\proof
Since $\DD(n)$ satisfies the localization property,
we can apply the argument of Thomason \cite[Prop. 2.8]{thomason} using
induction on the dimension of $X$, to reduce to showing
that for an artinian local ring $R$, essentially of finite type over $k$,
we have $\DD(n)(\Spec R)\cong R\Gamma(\Spec R_\et,\DD(n))$ for $n\leq 0$.
Since $\DD(n)(U)=\DD(n)(U^{red})$, we can assume that $R$ is
reduced, in which case it is the spectrum of a field $F$ of finite
transcendence degree $d$ over $k$. We have to show that the canonical map
$H_i(F,\DD(n))\to H_i(F_\et,\DD(n))$ is an
isomorphism for all $i$ and $n\leq 0$. Rationally,
Zariski and \'etale hypercohomology of the motivic complex agree.
With prime to $p$-coefficients, both sides agree for $i\geq d+n$
by the Beilinson-Lichtenbaum conjecture. For $i<d$ both sides
vanish because by Proposition \ref{identifyme},
$H_i(F_\et,\DD/m(n))=H^{2d-i}(F_\et,\mu_m^{\otimes (d-n)})$
and because the cohomological dimension of $F$ is $d$.
With mod $p$-coefficients, both sides agree for $n=0$ because
$\DD/p^r(0)\cong \nu_r^d[d]\cong R\epsilon_*\nu_r^d[d]\cong
R\epsilon_*\DD/p^r(0)$ from Proposition \ref{identifyme} and
\cite{katokuzumaki}, where $\epsilon$ denotes the map from
the small \'etale site to the Zariski site.
Finally, $\Z/p^r(n)=0$ for $n<0$.
\proofend

We remark that the same argument gives an unconditional result
for the cycle complex $\Z^c_{(p)}(n)$ localized at $p$.

\begin{corollary}\label{pushmap}
Let $f:X\to Y$ be a map over the perfect field $k$, and $n\leq 0$.
If $f$ is proper, then there is a functorial push-forward
$f_*:Rf_*\DD_X(n)\to \DD_Y(n)$ in the derived category of \'etale
sheaves. For arbitrary $f$, we obtain for every torsion sheaf $\f$
on $Y$ a functorial map
\begin{equation}\label{tracemap}
Rf_!(f^*\f\otimes \DD_X(n))\to \f\otimes \DD_Y(n).
\end{equation}
\end{corollary}

\proof
For proper $f$, the map is given as
$$Rf_*\DD_X(n) \stackrel{\sim}{\leftarrow} \epsilon^*Rf_*^\Zar \DD_X(n)^\Zar
\to \epsilon^*\DD_Y(n)^\Zar \cong \DD_Y(n).$$
The second map is the map induced by the proper push-forward of higher
Chow groups. The first map is the base-change map between push-forward on
the Zariski-site and push-forward on the \'etale site
$\epsilon^*f_*^{\Zar}\to f_*^{\et}\epsilon^*$.
The push-forward
on the Zariski site $\epsilon^*Rf_*^\Zar \DD_X(n)^\Zar$ and on the \'etale site
$Rf_*\DD_X(n)$ are the complexes of \'etale sheaves on $Y$ associated to the
complexes of presheaves
$U\mapsto R\Gamma_\Zar(f^{-1}U,\DD(n))$ and
$U\mapsto R\Gamma_\et(f^{-1}U,\DD(n))$, respectively.
Showing that these complexes are quasi-isomorphic is a problem local
for the \'etale topology, hence we can assume that $k$ is algebraically
closed. In this case, the base change map is a quasi-isomorphism by
Theorem \ref{locok}.

For arbitrary $f:X\to Y$, factor $f$ through
a compactification $X\stackrel{j}{\to} T\stackrel{g}{\to} Y$.
Writing $\f$ as a direct limit of $m$-torsion sheaves we can
assume that $\f$ is $m$-torsion and replace
$\f\otimes\DD_Y(n)$ by $\f\otimes_{\Z/m}\DD_Y/m(n)$. Then using
the proper base-change theorem, we obtain a map
\begin{multline*}
Rf_!(f^*\f\otimes_{\Z/m} \DD_X/m(n))=
Rg_*j_!j^*(g^*\f\otimes_{\Z/m}\DD_T/m(n))\\
\to Rg_*(g^*\f\otimes_{\Z/m} \DD_T/m(n))
\xleftarrow{\sim} \f\otimes^{\mathbb L}_{\Z/m}Rg_*\DD_T/m(n)\to
\f\otimes_{\Z/m} \DD_Y/m(n) .
\end{multline*}
The usual argument comparing compactifications
shows that the composition is independent of the compactification.
\proofend

Under the identification of Proposition \ref{locokp},
the trace map with mod $p^r$-coefficients and for $n=0$
was constructed by Jannsen-Saito-Sato \cite{JSS}.
If $f:X\to k$ is proper over a perfect field, then
for $n=0$, the trace map agrees on the stalk $\Spec \bar k$
with the map sending a complex to its
highest cohomology group, composed with the degree map,
$$tr: Rf_*\DD_X(\bar k)\stackrel{\sim}{\leftarrow}\DD_X(X_{\bar k})\to
CH_0(X_{\bar k})\stackrel{\deg_{\bar k}}{\longrightarrow}\Z.$$

\begin{corollary}\label{purity}
a) Let $i:Z\to X$ be a closed embedding with open complement $U$
over a perfect field $k$. Then for every $n\leq 0$,
we have a quasi-isomorphism
$\DD_Z(n)\xrightarrow{\sim} Ri^!\DD_X(n),$
or equivalently a distinguished triangle
$$\cdots \to  R\Gamma(Z_\et,\DD_Z(n))\to
R\Gamma(X_\et,\DD_X(n))\to R\Gamma(U_\et,\DD_U(n))\to\cdots  . $$

b) If $p:X\times \A^r\to X$ is the projection and $n\leq 0$, then we have
a quasi-isomorphism of complexes of \'etale sheaves
$Rp_*\DD_{X\times \A^r}(n)\cong \DD_X(n-r)[2r]$.
\end{corollary}

\proof
a) Over an algebraically closed field,
we get $\DD_Z(n)\xrightarrow{\sim} Ri^!\DD_X(n)$
because $\DD(n)(X)\cong R\Gamma(X_\et,\DD(n))$, and higher Chow groups satisfy
localization \cite{blocloc}.
The general case follows by applying $R\Gamma_{G}$ to the
distinguished triangle over the algebraic closure of $k$,
for $G$ the absolute Galois group of $k$.

b) This follows from the homotopy formula \eqref{affineform} for $\DD(n)$ by
Theorem \ref{locok}.
\proofend

As in Proposition \ref{zarnivlemma}, localization gives

\begin{corollary}\label{nivlemma}
If $X$ is a scheme over a perfect field $k$ and $n\leq 0$,
then there are spectral sequences
\begin{equation}\label{niveauet}
E^1_{s,t}=\bigoplus_{x\in X_{(s)}}H^{s-t}(k(x)_\et,\Z(s-n)) \Rightarrow
H_{s+t}(X_\et,{\B}(n)).
\end{equation}
In particular, $H_i(X_\et,{\B}(n))=0$ for $i<n$.
\end{corollary}

\subsection{An alternate argument}
We deduce a version of the previous results without
using the Beilinson-Lichtenbaum conjecture.
Instead, we use a theorem of Suslin \cite{suslinetale},
who proves that for a scheme $X$ of dimension $d$ over an algebraically closed
field $k$, and an integer $m$ not divisible by the characteristic
of $k$, there is an isomorphism of finite groups
$CH_0(X,i,\Z/m)\cong H^{2d-i}_c(X_\et,\Z/m)^*$. If $X$ is smooth,
then this implies by Poincar\'e duality
that the groups $H^j(X_\Zar,\Z/m(d))$ and $H^j(X_\et,\Z/m(d))$ are
abstractly isomorphic, but this is weaker then the Beilinson-Lichtenbaum
conjecture, which states that the canonical (change of topology) map is
an isomorphism.

\begin{proposition}\label{locc}
Let $k$ be a perfect field and $n\leq 0$.

a) For a closed embedding $i:Z\to X$, the canonical map
$\DD_Z(n)\to Ri^!\DD_X(n)$ is a quasi-isomorphism. In particular,
we obtain the spectral sequence \eqref{niveauet}.

b) If $p:X\times \A^r\to X$ is the projection, then we have
a quasi-isomorphism of complexes of \'etale sheaves
$Rp_*\DD_{X\times \A^r}(n)\cong \DD_X(n-r)[2r]$.

c) If $f:X\to k$ is proper, then we have a trace map
$Rf_*\DD_X\to \Z$ in the derived category of \'etale sheaves on $k$.
\end{proposition}

\proof
We first prove a) for $n=0$.
The statement is \'etale local, so we can assume that $k$ is
algebraically closed and that $Z$ and $X$ are strictly henselian.
Furthermore, it suffices to show the statement for smooth $X$.
Indeed, if we embed $j:X\to  T$ into a
smooth scheme, and if $\DD_X\to Rj^!\DD_T$ as well as the composition
$\DD_Z \to Ri^!\DD_X \to Ri^!Rj^!\DD_T$ are quasi-isomorphisms,
then so is the map $\DD_Z \to Ri^!\DD_X$.
By the remark after Theorem \ref{locok} and the proof of Corollary
\ref{purity}a), the statement is known
with $\DD_{(p)}$-coefficients, hence it suffices to prove it
for mod $m$-coefficients, with $m$ not divisible by the characteristic
of $k$.

Consider the following commutative diagram of maps of long
sequences, where the upper vertical
maps are the change of topology maps, and the map $g$ is induced
by the base-change $\epsilon^*Ri^!_\Zar \f\to Ri^!_\et\epsilon^*\f$
together with localization $Ri^!_\Zar\DD_X\cong \DD_Z$,
\begin{equation}\label{bigda}
\begin{CD}
@>>>H_i(Z_\Zar, \DD_Z/m)@>>> H_i(X_\Zar,\DD_X/m)@>>> H_i(U_\Zar, \DD_U/m)@>>>\\
@III @| @| @Vf_iVV \\
@>>>H_i(Z_\et,\DD_Z/m)@>>> H_i(X_\et,\DD_X/m)@>>> H_i(U_\et,\DD_U/m)@>>>\\
@III @Vg_iVV @| @| \\
@>>> H_i(Z_\et,Ri^!\DD_X/m)@>>> H_i(X_\et,\DD_X/m)@>>> H_i(U_\et,\DD_U/m)@>>>.
\end{CD}
\end{equation}
The left two upper maps are isomorphisms because $Z$ and $X$
are strictly henselian. The upper row is exact by localization
\cite{blocloc}, and the lower row is exact by definition.
By Proposition \ref{identifyme}, the groups $H_i(X_\Zar,\DD_X/m)$ vanish
for $i\not=2d$ because $X$ is strictly henselian. Hence for $i\not=2d, 2d-1$
a diagram chase in the commutative diagram
$$\begin{CD}
H_{i+1}(U_\Zar, \DD_U/m)@=H_i(Z_\Zar, \DD_Z/m)\\
@VfVV@|  \\
H_{i+1}(U_\et,\DD_U/m)@>>>H_i(Z_\et,\DD_Z/m)\\
@|@VgVV  \\
H_{i+1}(U_\et,\DD_U/m)@=H_i(Z_\et,Ri^!\DD_X/m)
\end{CD}$$
shows that all groups are isomorphic.
If $i=2d$, then by Proposition \ref{identifyme},
the map $H_{2d}(X_\et,\DD_X/m)\to H_{2d}(U_\et,\DD_U/m)$ is the
identity of the group $\Z/m$, and $H_{2d+1}(U_\et,\DD_U/m)=0$.
Consequently the lower row of \eqref{bigda} shows that
$H_{2d}(Z_\et,Ri^!\DD_X/m)=H_{2d-1}(Z_\et,Ri^!\DD_X/m)=0$.
The map $H_{2d}(Z_\et,\DD_Z/m)\to H_{2d}(X_\et,\DD_X/m)$
factors through $H_{2d}(Z_\et,Ri^!\DD_X/m)=0$, hence is the zero map.
Similarly, the map
$H_{2d+1}(U_\Zar, \DD_U/m)\to H_{2d}(Z_\Zar, \DD_Z/m)$ factors
through $H_{2d+1}(U_\et, \DD_U/m)=0$, hence is the zero map as
well and we can conclude that $H_{2d}(Z_\Zar, \DD_Z/m)=0$.
Finally, a part of the upper two rows of \eqref{bigda} gives
the commutative diagram with exact rows
$$\begin{CD}
\Z/m@>>> H_{2d}(U_\Zar, \DD_U/m)@>>>
H_{2d-1}(Z_\Zar, \DD_Z/m)@>>> 0\\
@|  @Vf_iVV \\
\Z/m@=H_{2d}(U_\et, \DD_U/m) \\
\end{CD}$$
We now invoke Suslin's theorem, which implies that the source and target of
the surjection $f_i$ have the same finite order to conclude that
$H_{2d-1}(Z_\Zar, \DD_Z/m)=0$.
This finishes the proof that the map $\Z^c_Z/m\to Ri^!\Z^c_X/m$
is a quasi-isomorphism.

We now prove b) for $n=0$.
Rationally and with $p$-primary coefficients, we have \'etale
hypercohomological descent, and the claim follows from \eqref{affineform}.
With prime to $p$-coefficients, we embed $X$ into a smooth scheme
and use localization a) to reduce to the case that $X$ is smooth.
Then we apply Proposition \ref{identifyme},
and use the homotopy invariance of \'etale cohomology of
$\mu_m^{\otimes (d+t-n)}$
(a consequence of the smooth base-change Theorem and the calculation
of \'etale cohomology of the affine line).

To obtain a) for arbitrary $n\leq 0$, we let $r=-n$, consider the diagram
$$\begin{CD}
Z\times \A^r@>i_1>> X\times \A^r\\
@Vp'VV @VpVV \\
Z@>i>> X
\end{CD}$$
and get that
\begin{multline*}
Ri^!\DD_X(n)\cong Ri^!Rp_*\DD_{X\times\A^r}(0)[2n]\\
\cong Rp'_*Ri_1^!\DD_{X\times \A^r}(0)[2n]
\cong Rp'_*\DD_{Z\times \A^r}(0)[2n]\cong \DD_Z(n).
\end{multline*}
Now b) for arbitrary $n\leq 0$ follows from this as above.

c) Since it suffices to find a map over the algebraic closure of $k$
compatible with the Galois action, we can suppose that $k$ is
algebraically closed; in this case $Rf_*\DD_X$ can be identified
with $R\Gamma(X_\et,\DD_X)$.
From localization a), we obtain a spectral sequence \eqref{niveauet},
$$E^1_{s,t}=\displaystyle\bigoplus_{x\in X_{(s)}}
H^{s-t}(k(x)_\et,\Z(s)) \Rightarrow H_{s+t}(X_\et,\DD_X).$$
The terms $E^1_{s,t}$ vanish for $t<0$ as in Theorem \ref{locok}
by reasons of cohomological dimension. Hence
$R^{-t}f_*\DD_X=H_t(X_\et,\DD_X)=0$ for $t<0$,
and the map $E^2_{0,0}\to R^0f_*\DD_X$ is an isomorphism.
On the other hand, comparing spectral sequences \eqref{niveau} and
\eqref{niveauet}, we obtain a diagram
$$\begin{CD}\displaystyle
\bigoplus_{x\in X_{(1)}}H^1(k(x),\Z(1))@>>>\displaystyle
\bigoplus_{x\in X_{(0)}}H^0(k(x),\Z(0))@>>> CH_0(X)\to 0\\
@VVV @VVV @VVV \\ \displaystyle
\bigoplus_{x\in X_{(1)}}H^1(k(x)_\et,\Z(1))@>>> \displaystyle
\bigoplus_{x\in X_{(0)}}H^0(k(x)_\et,\Z(0))@>>> E^2_{0,0}
\to 0.
\end{CD}$$
Since the left vertical maps are isomorphisms, so is the right vertical map.
We get the trace map as the map of complexes of Galois modules
$$tr:Rf_*\DD_X\stackrel{\sim}{\longleftarrow}
 \tau_{\leq 0}Rf_*\DD_X \to R^0f_*\DD_X
\stackrel{\sim}{\longleftarrow} E_{0,0}^2 \stackrel{\sim}{\longleftarrow}
CH_0(X) \stackrel{tr}{\longrightarrow} \Z .$$
\proofend

\section{The main theorem}

\begin{theorem}\label{fdual}
Let $f:X\to k$ be separated and of finite type over a perfect field.
Then for every constructible sheaf $\g$ on $X$,
there is a canonical quasi-isomorphism
$$R\Hom_X(\g,\DD_X)\cong  R\Hom_k(Rf_!\g,\Z).$$
\end{theorem}

\proof
We can assume that $k$ is algebraically closed, because
if $G$ is the Galois group of $\bar k/k$, then
$R\Hom_X(\g,\f)\cong R\Gamma_{G}R\Hom_{\bar X}(\g,\f)$.
Indeed, for an injective sheaf $\f$, $\Hom_{\bar X}(\g,\f)$
is flabby \cite[III Cor. 2.13c)]{milnebook}.
Note that $Rf_!\f =R\Gamma_c(X_\et,\f)$
over an algebraically closed field.
If $X\stackrel{j}{\to} T\stackrel{g}{\to} k$
is a compactification, then the trace map
of Proposition \ref{locc}c) induces by adjointness the pairing
\begin{multline*}\label{TYT}
\alpha(\g):R\Hom_X(\g,\DD_X)\cong R\Hom_T(j_!\g,\DD_T)\\
\to R\Hom_k(Rg_*j_!\g,Rg_*\DD_T)
\to R\Hom_k(Rf_!\g,\Z).
\end{multline*}
The standard argument comparing compactifications shows that
this does not depend on the compactification.
We proceed by induction on the dimension $d$ of $X$,
and assume that the theorem is known for schemes of dimension
less than $d$.

\begin{lemma}\label{L3}
If $f:U\to X$ is \'etale, then $\alpha(f_!\g)$ is a quasi-isomorphism on
$X$ if and only if $\alpha(\g)$ is a quasi-isomorphism on $U$.
In particular, if $f$ is finite and \'etale, then $\alpha(\g)$ is a
quasi-isomorphism if and only if $\alpha(f_*\g)$ is a quasi-isomorphism.
\end{lemma}

\proof
If $f:U\to X$ is \'etale, then since $Rf^!\DD_X=f^*\DD_X=\DD_U$,
$\alpha(f_!\g)$ can be identified with $\alpha(\g)$ by adjointness.
If $f$ is also finite, then $f_*=f_!$.
\proofend

\begin{lemma}\label{rrr}
Let $j:U\to X$ be a dense open subscheme of a scheme of dimension $d$.
Then $\alpha(\g)$ is a quasi-isomorphism if and only
if $\alpha(j^*\g)\cong \alpha(j_!j^*\g)$ is a quasi-isomorphism.
\end{lemma}

\proof
This follows by a $5$-Lemma argument and induction on the dimension
from the map of distinguished triangles
$$\begin{CD}
R\Hom_Z(i^*\g,\DD_Z)@>>>R\Hom_X(\g,\DD_X)
@>>> R\Hom_U(j^*\g,\DD_U)\\
@VVV @VVV @VVV\\
R\Hom(R\Gamma_c(Z_\et,i^*\g),\Z)@>>> R\Hom(R\Gamma_c(X_\et,\g),\Z)
@>>>R\Hom(R\Gamma_c(U_\et,j^*\g),\Z)\\
\end{CD}$$
arising by adjointness from the short exact sequence
$0\to j_!j^*\g\to \g\to i_*i^*\g\to 0$ and purity for $\DD_X$,
Proposition \ref{locc}a).
\proofend

\begin{lemma}\label{upper}
Let $X$ be a scheme of dimension $d$ over an algebraically closed
field $k$. Then for any constructible sheaf $\g$, we have 
$\Ext^i_X(\g,\Z^c)=0$ for $i>2d+1$.
\end{lemma}

\proof
If $\g$ is locally constant, then the stalk of $\EXT^i_X(\g,\Z^c)_x$
agrees with $\Ext^i_\Ab(\g_x,\Z^c_x)$by \cite[III Ex.1.31b]{milnebook},
hence vanishes for $i>1$ because $\Z^c$ is concentrated in non-positive
degrees. Since $X$ has cohomological dimension $2d$, we conclude
by the spectral sequence
\begin{equation}\label{ltgl}
E_2^{s,t}=H^s(X_\et,\EXT^t_X(\g,\f))\Rightarrow \Ext^{s+t}_X(\g,\f)
\end{equation}
in this case. In general, we proceed by induction on the dimension of $X$.
Let $j:U\to X$ be a dense open subset with complement $i:Z\to X$
such that $\g|_U$ is locally constant. The statement follows with the long
exact $\Ext$-sequence arising from the short exact sequence
$0\to j_!j^*\g \to \g \to i_*i^* \g\to 0$ by purity,
$$\cdots \to \Ext^i_Z(\g,\Z^c_Z(n))\to \Ext^i_X(\g,\Z^c_X(n))\to
\Ext^i_U(\g,\Z^c_U(n))\to \cdots.$$
\proofend

\begin{lemma}\label{L2}
If $\alpha(\g)$ is a quasi-isomorphism for all constant constructible
sheaves $\g$ on smooth and projective schemes,
then $\alpha(\g)$ is a quasi-isomorphism for all constructible
sheaves $\g$ on all schemes $X$.
\end{lemma}

\proof
Replacing $X$ by a compactification $j:X\to T$ and $\g$ by $j_!\g$, 
we can assume that $X$ is proper.
We fix $X$ and show by descending induction on $i$ that the map
$$ \alpha_i(\g):\Ext^i_X(\g,\DD_X)\to
R\Hom^i_k(R\Gamma(X,\g),\Z)
\cong H^{1-i}(X_\et,\g)^*$$
induced by $\alpha$ is an isomorphism for all constructible
sheaves $\g$ on $X$. By Lemma \ref{upper}, both sides vanish for large $i$.
We can find an alteration $\pi:Y\to X$ with $Y$ smooth and projective,
together with a dense open subset $j:U\to X$ such that the restriction
$p:V=\pi^{-1}U\to U$ is \'etale and $p^*j^*\g=C_V$ is constant:
$$\begin{CD}
V@>>> Y\\
@VpVV @V\pi VV \\
U@>j>> X.
\end{CD}$$
Indeed, find an open subset $U$ of $X$ and a finite \'etale cover
$V'$ of $U$ such that $\g|_{V'}$ is constant. Let $Y'$ be the closure
of $V'$ in $X\times T$ for some compactification $T$ of $V'$. Now let
$Y$ be a generically \'etale alteration of $Y'$ which is smooth and
projective, and shrink $U$ such that $V=V'\times_{Y'}Y=U\times_XY$
is finite and \'etale over $U$.

By hypothesis, $\alpha(C_Y)$ is a quasi-isomorphism, hence $\alpha(C_V)$ is
a quasi-isomorphism by
Lemma \ref{rrr} (note that $C_Y|_V=C_V$ because $Y$ is smooth),
and so is $\alpha(p_*C_V)$ by Lemma \ref{L3}.
By the proper base change theorem, $\pi_*\pi^*\g$ is constructible, and
$j^*\pi_*\pi^*\g=p_*p^*j^*\g=p_*C_V$. Thus
$\alpha(\pi_*\pi^*\g)$ is a quasi-isomorphism by Lemma \ref{rrr}.
Let $\g'$ be the cokernel of the adjoint inclusion $\g\to \pi_*\pi^*\g$, and
consider the map of long exact sequences
\begin{equation}\label{devis}
\begin{CD}
\Ext^i_X(\g',\DD_X)@>>>\Ext^i_X(\pi_*\pi^*\g,\DD_X)@>>>
\Ext^i_X(\g,\DD_X)\\
@V\alpha_i(\g') VV @| @V\alpha_i(\g)VV\\
H^{1-i}(X_\et,\g')^*@>>> H^{1-i}(X_\et,\pi_*\pi^*\g)^*
@>>> H^{1-i}(X_\et,\g)^*.
\end{CD}\end{equation}
If $\alpha_j(\g)$ is
an isomorphism for $j>i$ and all constructible sheaves $\g$,
then a $5$-Lemma argument in \eqref{devis} shows that $\alpha_i(\g)$
is surjective. Since this holds for all constructible sheaves,
in particular $\g'$, another application of the $5$-Lemma shows
that $\alpha_i(\g)$ is an isomorphism.
\proofend

It remains to prove

\begin{proposition}
Let $X$ be smooth and projective
over an algebraically closed field $k$. Then for every positive
integer $m$, the map $\alpha(\Z/m)$ is a quasi-isomorphism
$$ R\Hom_X(\Z/m,\DD_X)\xrightarrow{\sim}
R\Hom_k(Rf_*\Z/m,\Z).$$
\end{proposition}

\proof
We can assume that $X$ is irreducible of dimension $d$.
By Lemma \ref{intmodm}, the pairing agrees up to a shift with
$$\begin{CD}
R\Hom_{X,\Z/m}(\Z/m,\DD_X/m)@>>>
R\Hom_{k,\Z/m}(Rf_*\Z/m,\Z/m)\\
@| @| \\
R\Gamma(X_\et,\DD_X/m)@>>>
\Hom_\Ab(R\Gamma(X_\et,\Z/m),\Z/m).\\
\end{CD}$$
By induction on the number of prime factors of $m$ we can
assume that $m$ is a prime number.
If $m$ is prime to the characteristic of $k$, then we claim that
the lower row of the diagram agrees with Poincar\'e duality
\cite[VI Thm. 11.1]{milnebook}
$$
R\Gamma(X_\et,\mu_m^{\otimes d})[2d]\to\Hom_\Ab(R\Gamma(X_\et,\Z/m),\Z/m).
$$
Since $Rf_*\f=R\Hom_{X,\Z/m}(\Z/m,\f)$, it is easy to see that
our pairing agrees with the Yoneda pairing, and it remains to
show that the trace map agrees with the map of Proposition \ref{locc}c), i.e.
the following diagram commutes
$$\begin{CD}
R\Gamma(X_\et,\DD_X/m)@>f_*>> R\Gamma(k_\et,\DD_k/m)\\
@| @| \\
R\Gamma(X_\et,\mu_m^{\otimes d}[2d]) @> tr' >>  \Z/m.
\end{CD}$$
The trace map $tr'$ is characterized
by the property that it sends the class 
$i_*(1)\in H^{2d}(X_\et,\mu_m^{\otimes d})$ of a closed point $i:p\to X$
to $1$ \cite[VI Thm.11.1]{milnebook},
a property which is clear for the map of Proposition \ref{locc}c) by
functoriality. Hence it suffices to show the commutativity
of the diagram
$$\begin{CD}
\Z/m\cong R\Gamma(k_\et,\DD_k/m)@>i_* >>R\Gamma(X_\et,\DD_X/m)\\
@| @| \\
\Z/m\cong R\Gamma(k_\et,\Z/m)@>i_* >>
R\Gamma(X_\et,\mu_m^{\otimes d}[2d]),
\end{CD}$$
which follows from Proposition \ref{identifyme}.

If $m=p$ is the characteristic, then we claim that our pairing
agrees with
Milne's duality \cite{milnevalues}
\begin{equation}\label{dualp}
R\Gamma(X_\et,\nu^d)[d]\to  \Hom_\Ab(R\Gamma(X_\et,\Z/p),\Z/p).
\end{equation}
More precisely, consider the complexes of
commutative algebraic perfect $p$-torsion group schemes
$\underline{H}^\cdot(X,\nu^d)$ and $\underline{H}^\cdot(X,\Z/p)$ over $k$, see
\cite{milnevalues}. Finite generation of $H^i(X_\et,\Z/p)$ implies that the
unipotent part of $\underline{H}^i(X,\Z/p)$ is trivial, and the same then
holds for the unipotent part of $\underline{H}^i(X,\nu^d)$ by
loc. cit. Theorem 1.11. Hence the same Theorem shows that the Yoneda pairing
induces a duality of \'etale group schemes
$$ \underline{H}^\cdot(X,\nu^d)\to
\Hom(\underline{H}^\cdot(X,\Z/p),\Z/p)[-d],$$
and \eqref{dualp} are the global sections over $k$.
By the argument above, it suffices to show
that the following diagram commutes
$$\begin{CD}
R\Gamma(X_\et,\DD_X/p)@>f_*>> R\Gamma(k_\et,\DD_k/p)\\
@| @| \\
R\Gamma(X_\et,\nu^d[d]) @> tr' >>  \Z/p.
\end{CD}$$
Again the trace map is characterized
by the property that it sends the class of a closed point $i:p\to X$ to $1$
\cite[p.308]{milnevalues}, and it suffices to show the commutativity
of the diagram
$$\begin{CD}
\Z/p\cong H^0(k_\et,\DD_k/p)@>i_* >>H^0(X_\et,\DD_X/p)\\
@| @| \\
\Z/p\cong H^0(k_\et,\Z/p)@>i_* >>
H^d(X_\et,\nu^d),
\end{CD}$$
which again is Proposition \ref{identifyme}.
\proofend

\rem
If $X$ is smooth of dimension $d$ over a perfect field of
characteristic $p$, and $\g$ a locally constant $m$-torsion sheaf,
then we can identify the left hand side of Theorem \ref{fdual},
$$\Ext^i_X(\g,\DD_X)\cong H^{2d+i-1}(X_\et,\g^D),$$
where $\g^D =\HOM_X(\g,\mu_m^{\otimes d})$ if $p\not| m$,
and $\g^D=\HOM_X(\g,\nu_r^d)[d] $ if $m=p^r$.
Indeed, consider the the spectral sequence \eqref{ltgl}.
The calculation of the $\EXT$-groups is local for the \'etale
topology, and since $\g$ is locally constant, we can
calculate them at stalks \cite[III Ex.1.31b]{milnebook}
and assume that $\g=\Z/m$. Then by Lemma \ref{intmodm}
and Proposition \ref{identifyme},
$$\EXT^q_X(\Z/m,\DD_X) =
{\cal H}^{q-1}(\DD_X/m)\cong
\begin{cases}
\mu_m^{\otimes d} &p\not| m, q=1-2d \\
\nu_r^d &m=p^r, q=1-d\\
0 &\text{otherwise}.
\end{cases}$$

\begin{corollary}\label{adji}
Let $f:X\to Y$ be a map of schemes over a perfect field $k$
and $n\leq 0$.

a) For every locally constant constructible sheaf $\f$ on $Y$, the map
\eqref{tracemap} induces a quasi-isomorphism
$$ \DD_X(n)\otimes f^*\f\cong Rf^!(\DD_Y(n)\otimes\f) .$$

b) For every torsion sheaf $\g$ on $X$ and finitely generated,
locally constant sheaf $\f$ on $Y$, we have a functorial quasi-isomorphism
$$R\Hom_X(\g,\DD_X(n)\otimes f^*\f)\cong R\Hom_Y(Rf_!\g,\DD_Y(n)\otimes\f).$$
\end{corollary}

It was pointed out to us by A.Abbes that a) is false for general
constructible sheaves $\f$. For example, if $i$ is the inclusion of
a point into the affine line, $\f=i_*\Z/m$ and $n=0$, then the left hand
side is $\Z/m$, but the right hand side is
$Ri^!(\DD_{{\mathbb A}^1}\otimes i_*\Z/m)\cong
Ri^!(i_*\mu_m[2])\cong \mu_m[2]$ by Proposition \ref{identifyme}.

\medskip

\proof
a) Since the statement is local for the \'etale topology on $Y$, we can
assume that $\f$ is of the form $\Z/m$. If $Y=\Spec k$ and $n=0$, let
$j:U\to X$ be an \'etale map and $g=f\circ j$. Then it suffices to show
that the upper row in the following diagram is an isomorphism
$$\begin{CD}
R\Hom_U(\Z/m,\DD_U/m)@>ad(tr)>> R\Hom_U(\Z/m,j^*Rf^!\DD_k/m)\\
@VVV @| \\
R\Hom_k(Rg_!\Z/m,Rg_!\DD_X/m)@>tr>> R\Hom_k(Rg_!\Z/m,\DD_k/m).
\end{CD}$$
The diagram commutes by the property of adjoints, and the lower
left composition is a quasi-isomorphism by Theorem \ref{fdual}.
For arbitrary $Y$, if $t:Y\to k$ is the structure map, then by
the previous case, the second map and the composition in
$$ \DD_X/m \to Rf^!\DD_Y/m \to Rf^!Rt^!\DD_k/m $$
are quasi-isomorphisms, hence so is the first map.
If $n<0$, let $r=-n$ and consider the following
commutative diagram
$$\begin{CD}
X\times{\Bbb A}^r@>p'>> X\\
@VgVV @VfVV \\
Y\times {\Bbb A}^r @>p>> Y.
\end{CD}$$
By homotopy invariance $Rp_*\DD_{Y\times{\Bbb A}^r}(0)[2n]\cong \DD_Y(n)$
and similar for $X$, and we obtain
\begin{multline*}
\DD_X(n)\cong Rp'_*\DD_{X\times{\Bbb A}^r}[2n]
\cong Rp'_*Rg^!\DD_{Y\times{\Bbb A}^r}[2n]\\
\cong Rf^!Rp_*\DD_{Y\times{\Bbb A}^r}[2n]\cong Rf^!\DD_Y(n).
\end{multline*}

b) If $\g$ is a constructible $m$-torsion sheaf, we can
replace $\f$ by $\f/m$ and use a).
For general $\g$, write $\g=\colim \g_i$ as a filtered colimit of
constructible sheaves \cite[II Prop.0.9]{adt}. Then the canonical map
$$ R\lim R\Hom_X(\g_i,\DD_X(n)\otimes f^*\f)\to
R\lim R\Hom_Y(Rf_!\g_i,\DD_Y(n)\otimes \f)$$
induces a map of spectral sequences \cite{roos}
$$ \begin{CD}
E_2^{s,t}=\lim^s \Ext^t_X(\g_i,\DD_X(n)\otimes f^*\f) @EEE \Rightarrow
\Ext^{s+t}_X(\colim \g_i,\DD_X(n)\otimes f^*\f) \\
@VVV @VVV \\
E_2^{s,t}=\lim^s \Ext^t_Y(Rf_!\g_i,\DD_Y(n)\otimes \f) @EEE \Rightarrow
\Ext^{s+t}_Y(\colim Rf_!\g_i,\DD_Y(n)\otimes \f) .
\end{CD}$$
The map on $E_2$-terms is an isomorphism by the above,
and by the following Lemma
the spectral sequences converge. Hence the map on abutments
is an isomorphism. Finally, \'etale cohomology with compact support
commutes with filtered colimits.
\proofend

\begin{lemma}\label{lower}
Let $\g$ be a torsion sheaf and $\f$ be finitely generated locally constant
sheaf on a scheme $X$ of dimension $d$ over an algebraically closed field.
Then $\Ext^i_X(\g,\Z^c_X(n)\otimes \f)=0$ for $i<-2d$ and $n\leq 0$.
\end{lemma}

\proof
By \eqref{ltgl} it suffices to show that $\EXT^i_X(\g,\Z^c_X(n)\otimes \f)=0$
for $i<-2d$,
and since this is \'etale local, we can assume that $\f=\Z$ or $\f=\Z/m$.
By the long exact coefficient sequence, both cases follows from
$\EXT^i_X(\g,\Z^c_X(n))=0$ for $i<-2d+1$. Since $\g$ is torsion, this will
follow if $\EXT^i_X(\g,\Q/\Z^c_X(n))=0$ for $i<-2d$, which by
\cite[III Rem.1.24]{milnebook} will follow from 
$\Ext^i_V(\g|_V,\Q/\Z^c_V(n))=0$ for all $V$ \'etale over $X$.
Let $U$ be a dense smooth open subscheme of $V$ with complement
$Z$. Then we conclude by induction on the dimension of $V$,
Proposition \ref{identifyme} and the long exact localization sequence
arising from purity
$$\to \Ext^i_Z(\g|_Z,\Q/\Z^c_Z(n))\to \Ext^i_V(\g|_V,\Q/\Z^c_V(n))\to
\Ext^i_U(\g|_U,\Q/\Z^c_U(n)) \to .$$
\proofend

\begin{corollary}\label{betterupper}
Let $\g$ be a torsion sheaf on a scheme $X$ over a perfect field $k$,
and let $n\leq 0$. Then for any constructible sheaf $\g$,
$\Ext^i_X(\g,\Z^c_X(n))$ and $\EXT^i_X(\g,\Z^c_X(n))$ vanish
for $i>\cd k+1$.
In particular, $\Z^c(n)$ has quasi-injective dimension $\cd k+1$
in the sense of SGA 5 I Def.1.4.
\end{corollary}

\proof
If $k$ is algebraically closed, then
$$\Ext^i_X(\g,\Z^c_X(n))\cong \Ext^i_\Ab(R\Gamma_c(X_\et,\g),\Z^c_k(n))
\cong \Hom(H^{1-i}_c(X_\et,\g),\Q/\Z^c_k(n)).$$
By Proposition \ref{identifyme}, the complex $\Q/\Z^c_k(n)$
is concentrated in degree zero, hence this vanishes for $i>1$.
In the general case, we use the spectral sequence
$$ H^s(Gal(k), \Ext^t_{\bar X}(\g,\Z^c_X(n)))\Rightarrow
\Ext^{s+t}_X(\g,\Z^c_X(n)).$$
The statement for the extension sheaves follows because
$\EXT^i_X(\g,\Z^c(n))$ is the sheaf associated to the presheaf
$U\mapsto \Ext^i_U(\g|_U,\Z^c_U(n))$.
\proofend

\section{Duality}
If a perfect field $k$ has duality for Galois cohomology with a dualizing
sheaf that is related to some $\DD_k(n)$, then Theorem \ref{fdual}
gives a duality theorem over $k$. For example,
if $k$ is algebraically closed, then we immediately obtain from
Theorem \ref{fdual} a quasi-isomorphism
$$R\Hom_X(\g,\DD_X)\cong  R\Hom_\Ab(R\Gamma_c(X_\et,\g),\Z)$$
for every constructible sheaf $\g$ on $X$. In particular, we get perfect
pairings of finitely generated groups
$$\Ext^{1-i}_X(\g,\DD_X)\times H^i_c(X_\et,\g)\to \Q/\Z.$$
From Lemma \ref{intmodm} and Theorem \ref{locok} this gives
\begin{equation}\label{uyh}
CH_0(X,i,\Z/m)\cong H^i_c(X_\et,\Z/m)^*,
\end{equation}
generalizing Suslin's theorem \cite{suslinetale} to arbitrary $m$.

\subsection{Finite fields}
Let $G$ be the absolute Galois group of $\bar\F_q/\F_q$ with $q=p^r$.
The following theorem has been proved by Jannsen-Saito-Sato \cite{JSS}
and Moser \cite{moser} for $p$-power torsion sheaves, by Deninger
\cite[Thms. 1.4, 2.3]{deninger2} for curves, and smooth
schemes over curves and coefficients of order prime to $p$,
and by Spie\ss\  \cite{spiess} for surfaces.

\begin{theorem}\label{finitedual}
Let $X$ be a scheme over a finite field, and
$\g$ be a torsion sheaf on $X$. Then
there is a quasi-isomorphism
$$R\Hom_X(\g,\DD_X)\xrightarrow{\sim}R\Hom_\Ab(R\Gamma_c(X_\et,\g),\Z)[-1].$$
In particular, if $\g$ is constructible, there are perfect pairings
of finite groups
$$\Ext^{2-i}_X(\g,\DD_X)\times H^i_c(X_\et,\g)\to \Q/\Z.$$
\end{theorem}

\proof
If $\bar f:X\times_{\F_q}\bar \F_q\to \bar \F_q$ is the structure map,
and if we apply $R\Gamma_{G}$ to the pairing over $\bar\F_q$,
then we get
$$ R\Hom_X(\g,\DD_{\bar X})=
R\Gamma_{G}R\Hom_{\bar X}(\g,\DD_{\bar X})
\cong R\Gamma_{G}R\Hom_{\bar \F_q}(R\bar f_!\g,\Z).$$
By duality for Galois cohomology, this is quasi-isomorphic to
$$ R\Hom_\Ab(R\Gamma_{G}R\bar f_!\g,\Z)[-1]
\cong R\Hom_\Ab(R\Gamma_c(X_\et,\g),\Z)[-1].$$
\proofend

Let $H_i^K(X,\Z/m)$ be the $i$th homology group of the Kato complex \cite{kato}
\begin{equation}\label{katocomplex}
\oplus_{X_{(0)}} H^{1}(k(x)_\et,\Z/m) \leftarrow
\oplus_{X_{(1)}} H^{2}(k(x)_\et,\Z/m(1))\leftarrow \cdots.
\end{equation}
Kato \cite[Conj.0.3]{kato} conjectures that $H_i^K(X,\Z/m)=0$
for $i>0$ and $X$ smooth and proper.
Kato's conjecture has been proved in low degrees by Colliot-Th\'el\`ene
\cite{colliet}, and in general by Jannsen and Saito
\cite{jannsen, jannsensaito} assuming resolution of singularities.
One important application of Kato homology is, in view of
$H^1(X_\et,\Z/m)^*\cong \pi_1^{ab}(X)/m$, the following

\begin{corollary}
Assuming the Beilinson-Lichtenbaum conjecture, there is,
for every scheme over a finite field, an exact sequence
\begin{multline}\label{S}
\cdots \to CH_0(X,i,\Z/m)\to
H^{i+1}_c(X_\et,\Z/m)^* \to H_{i+1}^K(X,\Z/m)\to \cdots.
\end{multline}
\end{corollary}

\proof
By Theorem \ref{finitedual} and Lemma \ref{intmodm}, we have
\begin{multline*}
H^{i+1}_c(X,\Z/m)^*=\Ext^{1-i}_X(\Z/m,\DD_X)\\ \cong
\Ext^{-i}_{X,\Z/m}(\Z/m,\DD_X/m)\cong H_i(X_\et,\DD_X/m).
\end{multline*}
The Corollary follows by comparing the niveau spectral sequences
\eqref{niveau} and \eqref{niveauet} with $\Z/m$-coefficients,
and identifying the terms $E_{s,t}^1$ with $t\geq n$ by
the Beilinson-Lichtenbaum conjecture.
\proofend

\subsection{Local fields}
Let $k$ be the field of fractions of a henselian discrete valuation ring
of characteristic $0$ with finite residue field, for example a local field
of mixed characteristic, and $G$ the Galois group of $k$.

\begin{theorem}\label{localdual}
If $f:X\to k$ is a scheme over $k$, and $\g$ a torsion sheaf, then
there is a quasi-isomorphism
$$ R\Hom_X(\g,\DD_X(-1))
\cong R\Hom_{Ab}(R\Gamma_c(X,\g),\Z)[-2].$$
In particular, for constructible $\g$, we have perfect pairings
of finite groups
$$\Ext^{3-i}_X(\g,\DD_X(-1))\times H^i_c(X_\et,\g)\to \Q/\Z.$$
\end{theorem}

\proof
From Corollary \ref{adji} we get the following quasi-isomorphisms
\begin{multline*}
R\Hom_X(\g,\DD_X(-1))\cong
R\Gamma_{G}\Hom_{\bar X}(\g, \DD_{\bar X}(-1))\\
\cong R\Gamma_{G}R\Hom_{\bar k}(R\Gamma_c(\bar X,\g),\DD_{\bar k}(-1))
\cong R\Hom_{G}(R\Gamma_c(\bar X,\g),\DD_{\bar k}(-1)).
\end{multline*}
The claim follows with the following Lemma.
\proofend

\begin{lemma}(Duality for Galois cohomology)
If $G$ is the Galois group of $\bar k/k$, and $C^\cdot$
a complex of torsion $G$-modules, then there is
a quasi-isomorphism
$$ R\Hom_{G}(C^\cdot,\DD_k(-1))\cong
R\Hom_\Ab(R\Gamma(G,C^\cdot),\Z)[-2].$$
\end{lemma}

\proof
Since $\DD_k(-1)\cong {\Bbb G}_m[-1]$, this is local duality
$$ R\Hom_{G}(C^\cdot,\G_m)[-1]\cong R\Hom_\Ab(R\Gamma(G,C^\cdot),\Q/\Z)[-3].$$
Indeed, the Yoneda pairing induces an isomorphism \cite[I Thm. 2.14]{adt}
$$\Ext^i_{G}(M,\G_m)\cong \Hom_\Ab(H^{2-i}(G,M),\Q/\Z)$$
for every $G$-module $M$, and the result for complexes follows
by a spectral sequence argument.
\proofend

\rem
If $X$ is smooth of dimension $d$ over a local field $k$
of characteristic $0$, $\g$ is a locally constant $m$-torsion sheaf,
and $\g^D:= \HOM(\g,\mu_m^{d+1})$, then we recover the perfect pairing
of finite groups
$$H^{2d+2-i}(X_\et,\g^D)\times H^i_c(X_\et,\g)\to \Q/\Z.$$
This follows from Theorem \ref{localdual}, using Proposition \ref{identifyme},
and the degeneration of the local-to-global spectral sequence of
Ext-groups.

Let $H_i^K(X,\Z/m)$ be the $i$th homology of the Kato complex
\begin{equation}\label{katolocal}
\oplus_{X_{(0)}} H^2(k(x)_\et,\Z/m(1)) \leftarrow
\oplus_{X_{(1)}} H^3(k(x)_\et,\Z/m(2))\leftarrow \cdots.
\end{equation}
Kato conjectures \cite[Conj. 5.1]{kato} that the terms
$H_i^K(X,\Z/m)$ vanish for $X$ smooth and proper and $i>0$.
The same proof as for \eqref{S} gives

\begin{corollary}
Assuming the Beilinson-Lichtenbaum conjecture,
there is, for every scheme over a local field of characteristic $0$,
a long exact sequence,
\begin{equation}\label{T}
\cdots \to CH_{-1}(X,i,\Z/m)
\to H^i_c(X_\et,\Z/m)^* \to H_i^K(X,\Z/m)\to \cdots.
\end{equation}
\end{corollary}

\rem If $K$ is a $r$-local field of characteristic $0$ such that
the finite field in the definition of $K$ has characteristic $p$,
then one can use local duality for its Galois cohomology to prove
a quasi-isomorphism
$$ R\Hom_X(\g,\DD_X(-r))
\cong R\Hom_{Ab}(R\Gamma_c(X,\g),\Z)[-r-1]$$
for torsion sheaves $\g$ of order not divisible by $p$,
generalizing a result of Deninger-Wingberg \cite{deningerwingberg}.

\section{Applications}
\subsection{Rojtman's theorem}
We generalize Rojtman's theorem  \cite{milneroitman, roitman}
to normal projective varieties. Our proof follows the line of Bloch
\cite{blochrojt} and Milne \cite{milneroitman}.
Let $X$ be a proper scheme over an algebraically closed field.
Then by \cite[Prop. 4.16]{milnebook} there is an isomorphism
$H^1(X_\et,\Z/m)\cong \Hom_k(\mu_m,\Pic_X)$, where $\Hom_k$ is the group of
homomorphisms of flat group schemes over $k$.
Let $\Pic^\tau_X$ be the scheme representing line  bundles,
such that a power is algebraically equivalent to $0$. Then the quotient
$\Pic_X/\Pic^\tau_X$ is a finitely generated free group,
hence $\Hom_k(\mu_m,\Pic_X)\cong \Hom_k(\mu_m,\Pic^\tau_X)$.
The quotient $C=\Pic^\tau_X/\Pic^{0,\red}_X$ by the reduced part of
the connected component is a finite
group scheme, and $ \Pic^{0,\red}_X$ is a smooth commutative
algebraic group, which by Chevalley's theorem is an extension
of an abelian variety by a linear algebraic group.
Since $\Ext_k(\mu_m,G)=0$ for every smooth connected group scheme,
we obtain a short exact sequence
$$0\to \Hom_k(\mu_m,\Pic^{0,\red}_X)\to \Hom_k(\mu_m,\Pic_X)
\to \Hom_k(\mu_m, C)  \to 0.$$
If $(_m\Pic^{0,\red}_X)^\vee$ is the Cartier dual of the $m$-torsion part
${}_m\Pic^{0,\red}_X$,
then by Cartier-Nishi duality, we obtain a short exact sequence
$$0\to \Hom_k((_m\Pic^{0,\red}_X)^\vee,\Z/m)\to \Hom_k(\mu_m,\Pic_X)
\to  \Hom_k(C^\vee,\Z/m)  \to 0.$$
Since $\Z/m$ is \'etale, $\Hom_k(G,\Z/m)\cong \Hom_\Ab(G(k),\Z/m)$ for every
finite group scheme $G$, hence taking Pontrjagin duals we arrive at
\begin{equation}\label{neron}
0\to  C^\vee(k)/m \to H^1(X_\et,\Z/m)^*
\to ({}_m\Pic^{0,\red}_X)^\vee(k) \to 0.
\end{equation}
Let $\Alb_X$ be the Albanese variety
in the sense of Serre \cite{serre}, i.e. the universal
object for maps from $X$ to abelian varieties.
The universal map $X\to \Alb_X$ induces the albanese map
$CH_0(X)^0 \to \Alb_X(k)$, which is covariantly functorial for maps
between normal schemes.

\begin{theorem}\label{rojtman}
Let $X$ be a normal scheme, projective over an algebraically
closed field. Then the albanese map induces an isomorphism
$$ {}_\tor CH_0(X)\xrightarrow{\sim}{}_\tor \Alb_X(k), $$
and $CH_0(X,1)\otimes\Q/\Z=0$.
\end{theorem}

The theorem re-proves and generalizes (to include the $p$-part) a result of
S.Saito \cite{shuji}.
Our result differs from the results of Levine and Krishna-Srinivas
\cite{levinerojtman, srin},
who compare the torsion in the (cohomological) Chow group of
Levine-Weibel \cite{levineweibel} to the torsion of the Albanese variety
in the sense of Lang-Weil, i.e. the universal object for {\it rational}
maps from $X$ to abelian varieties.

\smallskip

\proof
For normal, projective $X$, $\Pic^{0,\red}_X$ is an abelian variety by
\cite[Cor.3.2]{fga} or \cite[Rem. 5.6]{kleiman},
and its dual $(\Pic^{0,\red}_X)^t$ is $\Alb_X$
by \cite[Thm.3.3]{fga} or \cite[Rem. 5.25]{kleiman}.
In particular, $(_m\Pic^{0,\red}_X)^\vee={}_m\Alb_X$.
The usual argument of induction on the dimension
\cite[Lemma 2.1]{milneroitman} shows that the albanese map is surjective
on $l^n$-torsion for every $l$,
because a generic hyperplane section of a normal scheme is again
normal \cite{seidenberg}. Consider the diagram with exact rows
arising from \eqref{uyh} and \eqref{neron},
$$\begin{CD}
CH_0(X,1)\otimes \Z/l^n@>>> CH_0(X,1,\Z/l^n)@>>> {}_{l^n}  CH_0(X)@>>> 0\\
@III @| \\
C^\vee(k)/l^n @>>> H^1(X_\et,\Z/l^n)^* @>>>  {}_{l^n}\Alb_X(k)@>>> 0.
\end{CD}$$
Since $C$ is a finite group scheme, the lower left
term vanishes in the colimit. Counting (finite) coranks we obtain
$CH_0(X,1)\otimes\Q_l/\Z_l=0$, and an isomorphism
${}_\tor CH_0(X)\cong {}_{\tor}\Alb_X(k)$. Since the albanese map
is a surjection of divisible $l$-torsion groups of the same
corank, it must be an isomorphism.
\proofend

If $X$ is an arbitrary proper scheme, we define
$M^1(X)$ to be the quotient of $\Pic^{0,\red}_X$ by its largest
unipotent subgroup. Then we have a short exact sequence
\begin{equation}\label{1mot}
0\to T_X \to M^1(X)\to A_X \to 0,
\end{equation}
where $T_X$ is a torus and $A_X$ an abelian variety.
We let $M_1(X)$ be the dual of $M^1(X)$ in the category of $1$-motives;
see \cite[\S 10]{delignehodgeIII} for the definition and basic properties.
Concretely, $M_1(X)$
is the complex $\chi(T_X)\stackrel{f}{\to} A_X^t(k)$, where $\chi(T_X)$ is the
character group of $T_X$, $A^t_X$ the dual abelian variety of $A_X$,
and $f$ is given by pushing-out \eqref{1mot} along a given map
$T_X\to {\Bbb G}_m$ to obtain an element in $\Ext^1(A_X,{\Bbb G}_m)=A_X^t(k)$.
The Tate realization $T_m(M)$ of a $1$-motive $M=[F\to G]$
is the cone of multiplication by $m$. Since abelian varieties and
tori are divisible, $T_m(M)$ is concentrated in a single degree, and there is
a short exact sequence
$ 0\to {}_m G(k)\to T_m(M)\to F/m \to 0$.
In our case, we obtain a short exact sequence
\begin{equation}\label{albses}
0\to {}_\tor A_X^t(k)\to \colim_m T_mM_1(X)\to \chi(T_X)\otimes\Q/\Z
\to 0.
\end{equation}
If $M$ and $M^\vee $ are dual $1$-motives, then Cartier duality gives a
perfect pairing $T_m(M)\times T_m(M^\vee)\to \mu_m$.

\begin{proposition}
Let $X$ be a proper scheme over an algebraically closed field.
Then there is an isomorphism
$$CH_0(X,1,\Q/\Z)\cong \colim T_mM_1(X).$$
\end{proposition}

\proof
There are no homomorphisms and extensions
between $\mu_m$ and a unipotent group, hence
$$\Hom_k(\mu_m,\Pic^{0,\red}_X)\cong \Hom_k(\mu_m,M^1(X))
\cong \Hom_k(\mu_m,T_mM^1(X)),$$
and the argument which leads to \eqref{neron} gives
a short exact sequence
$$0\to  C^\vee(k)/m \to CH_0(X,1,\Z/m) \to T_mM_1(X) \to 0.$$
The result follows because, in the colimit, the first term vanishes.
\proofend

\example
Let $C=E/0\sim P$ be an elliptic curve $E$ with the points
$0$ and $P$ identified to the point $Q\in C$.
Then from the long exact sequences for higher Chow groups and \'etale
cohomology arising from the blow-up diagram
$$\begin{CD}
\{P,0\}@>>> E\\
@VVV @VVV\\
Q @>>> C
\end{CD}$$
we obtain $H^1(C_\et,\Z)=\Z$, hence $\chi(T_C)\cong \Z$ \cite{ichpicard},
and $CH_0(C,1,\Q_l/\Z_l)$ has corank $3$ for $l\not=\chr k$.
If $P$ is not a torsion point, then
$CH_0(C,1)\otimes\Q_l/\Z_l=0$ and ${}_{l^\infty}CH_0(C)$ has
corank $3$. Hence the maps
\begin{align*}
{}_\tor A_C^t(k) &\xrightarrow{f} {}_\tor CH_0(C)\\
CH_0(C,1)\otimes\Q/\Z&\xrightarrow{g} \chi(T_C)\otimes\Q/\Z ,
\end{align*}
obtained from \eqref{albses} have infinite cokernels,
and ${}_\tor CH_0(C)$ cannot be parametrized by an abelian
variety because it has odd corank.
If $P$ is an $m$-torsion point, then $CH_0(C,1)\otimes\Q/\Z=\Q/\Z$
and $f$ and $g$ are surjective with kernel $\Z/m$.

\subsection{The fundamental group}
Let $X$ be a proper scheme over an algebraically
closed field and consider the (profinite) fundamental group.
Then we have an isomorphism
$$ \pi_1^{ab}(X)=H^1(X_\et,\Q/\Z)^*\cong CH_0(X,1,\hat \Z) ,$$
in particular, a short exact sequence
\begin{equation}\label{roit}
0\to CH_0(X,1)^\wedge \to \pi_1^{ab}(X)
\to TCH_0(X)\to 0.
\end{equation}
If $X$ is normal, then taking the inverse limit in \eqref{neron},
and using Thereom \ref{rojtman}, we see that this sequence agrees
with the sequence
\begin{equation}\label{klang}
0\to  C^\vee(k) \to \pi_1^{ab}(X) \to T\Alb_X(k)\to 0,
\end{equation}
in particular, $CH_0(X,1)^\wedge$ is finite. The latter sequence
can be found in Milne \cite[Cor. III 4.19]{milnebook}, and in
Katz-Lang \cite[Lemma 5]{katzlang}.

For geometrically connected $X$ over a perfect field $k$,
the dual of the Hochschild-Serre spectral sequence 
$H^s(k,H^t(\bar X_\et,\Q/\Z))\Rightarrow H^{s+t}(X_\et,\Q/\Z)$
gives an exact sequence
\begin{equation}\label{despi}
H^2(k_\et,\Q/\Z)^*\to
\pi_1^{ab}(X\times_k\bar k)_{\Gal(k)}\to
\pi_1^{ab}(X)\to \Gal(k)^{ab} \to 0 ,
\end{equation}
which is short exact if $X$ has a $k$-rational point
\cite[Lemma 1]{katzlang}, or if
$k$ is finite, local, or global, because then $H^2(k_\et,\Q/\Z)=0$.
By the argument in \cite{katzlang}, the sequence \eqref{klang} implies that
$\pi_1^{ab}(X\times_k\bar k)_{\Gal(k)}$ is finite if
$k$ is absolutely finitely generated and $X$ is proper and normal.

Let $X$ be a proper scheme over a  finite field $k$ with Galois group $G$.
Let $\pi_1^{ab}(X)^0$ be the kernel of $\pi_1^{ab}(X)\to G$,
and $CH_0(X)^0$ the subgroup of cycles of degree zero of $CH_0(X)$.

\begin{proposition}
If $X$ is proper and geometrically connected over a finite field,
then we have a short exact sequence
$$0\to CH_0(\bar X,1)^\wedge_{G}\to\pi_1^{ab}(X)^0\to
(CH_0(\bar X)^0)^G\to 0.$$
If $X$ is normal, then these groups are finite, and the right hand
term agrees with $\Alb_X(k)$.
\end{proposition}

\proof
By \eqref{despi}, it suffices
to calculate the cokernel of $1-F$ on the sequence \eqref{roit},
for $F$ the Frobeinus, and to show that 
$TCH_0(\bar X)_{1-F}\cong CH_0(\bar X)^G$.
We can replace $CH_0(\bar X)$ by the divisible group $CH_0(\bar X)^0$
in \eqref{roit}.
Since the finitely generated group $CH_0(X)^0$ surjects onto
$(CH_0(\bar X)^0)^G$, the Galois-coinvariants of $CH_0(\bar X)^0$ are divisible
and finite, hence vanish. The short exact sequence
$$0\to (CH_0(\bar X)^0)^G\to CH_0(\bar X)^0
\stackrel{1-F}\longrightarrow CH_0(\bar X)^0\to 0$$
gives rise to the exact sequence
$$0\to {}_m(CH_0(\bar X)^0)^{G}\to {}_mCH_0(\bar X)^0
\stackrel{1-F}\longrightarrow  {}_mCH_0(\bar X)^0\to
(CH_0(\bar X)^0)^G/m \to 0.$$
The result follows by taking limits.
\proofend

\subsection{Duality theory}
For a complex $\f^\cdot$ of torsion sheaves on a scheme $X$ over a 
perfect field $k$, we consider the functor
$$\du(\f^\cdot )=R\HOM_X(\f^\cdot,\DD_X)$$
(if $\f^\cdot$ is unbounded, this is defined using a
K-injective resolution as in \cite[Prop.6.1]{spaltenstein}).

\begin{proposition}
The functor $\du$ sends bounded above complexes to bounded
below complexes and conversely.
\end{proposition}

\proof
Since the statement is \'etale local, we can assume that $k$ is
algebraically closed. By Lemma \ref{lower} and Corollary \ref{betterupper},
$\EXT^i_X(\g,\DD_X)$ vanishes for any torsion 
sheaf $\g$ unless $-2d\leq i\leq 1$. Hence the spectral sequence
$$ E_2^{p,q}=\EXT^p_X({\cal H}^{-q}(\f^\cdot),\DD_X)
\Rightarrow \EXT^{p+q}_X(\f^\cdot,\DD_X),$$
converges, and the claim follows.
\proofend

\begin{theorem} (Exchange formulas)\label{exchange1}
Let $f:X\to Y$ be a map between schemes over a perfect field
and $\f$ and $\g$ be constructible sheaves on $X$.
Then the following formulas hold
\begin{align*}
\du(\g\otimes \f)&\cong R\HOM(\g,\du(\f))\\
Rf_*\du(\f) &\cong \du (Rf_!\f) \\
Rf^!\du(\f) &\cong \du (f^*\f).
\end{align*}
\end{theorem}

\proof
The first formula is adjointness of $\Hom$ and $\otimes$.
For the second formula, in view of
$$R\HOM(\g,\f)(U)\cong R\Hom_U(\g|_U,\f|_U)\cong R\Hom_X(j_!j^*\g,\f),$$
it suffices to prove
$ R\Hom_X(\g,\DD_X)\cong R\Hom_Y(Rf_!\g,\DD_Y)$, for
constructible $\g$, which is Corollary \ref{adji}.
The last formula holds by SGA 4 XVIII Cor. 3.1.12.2.
\proofend

\begin{theorem} (Biduality)
Let $\f$ be a constructible $m$-torsion sheaf for some integer $m$ not
divisible by the characteristic of $k$. Then
$$ \f\cong \du(\du(\f)).$$
In particular,
\begin{align*}
\du (\g\otimes \du(\f)) &\cong R\HOM(\g,\f) \\
Rf_!\du(\f) &\cong \du (Rf_*\f) \\
f^*\du(\f) &\cong \du(Rf^!\f).
\end{align*}
\end{theorem}

\proof
For a closed embedding $i:Z\to X$ with open complement $j:U\to X$,
$R\HOM(j_!\g,\f)\cong j_!R\HOM(\g,j^*\f)$,
hence $\du(\du(j_!j^*\f))\cong  j_!\du(\du(j^*\f))$
on the one hand, and $R\HOM(i_*\g,\f)\cong i_*R\HOM(\g,Ri^!\f)$
hence $\du(\du(i_*i^*\f))\cong  i_*\du(\du(i^*\f))$
on the other hand. Thus we can use devissage to reduce to the case
that $\f\cong \Z/m$ on a smooth and proper scheme $X$
of dimension $d$ over $k$. In this case, by Lemma \ref{intmodm},
the statement reduces to
$ \Z/m\cong R\HOM_{\Z/m}(\mu_m^{\otimes d},\mu_m^{\otimes d})$.
Since $\mu_m$ is locally constant, this can be
checked at stalks, where it is clear.
The other formulas follow from this by
substitution in Theorem \ref{exchange1}, see SGA 5 I Prop.1.12.
\proofend

\rem
The biduality map at the characteristic is not an isomorphism in general.
For example, if $k$ is algebraically closed and $U$ is an open subset
of the projective line $\P^1_k$, then the localization
sequence for cohomology with compact support gives a short exact
sequence
$$ 0\to \bigoplus_{x\in \P^1-U}i^*_x\nu^1 \to H^1_c(U_\et,\nu^1)\to \Z/p\to 0$$
and it follows that the dual
$\HOM_{\Z/p}(\nu^1,\nu^1)(U)=
\Hom_{\P^1,\Z/p}(j_!\nu^1,\nu^1)$ of $H^1_c(U_\et,\nu^1)$
is very large, hence $\Z/p\not= R\HOM_{\Z/p}(\nu^1,\nu^1)$.
The small \'etale site does not give well-behaved
extension groups for non-constructible sheaves.
See Kato \cite{katoduality} for a good duality using the relative perfect
site.

\section{One dimensional bases}
For the remainder of the paper we assume the validity of the
Beilinson-Lichtenbaum conjecture.
Let $S$ be a connected one-dimensional regular scheme
such that all closed points have perfect residue fields, and the
generic point $\eta$ has characteristic $0$ (the case that $S$ is a curve
over a finite field is covered by Theorem \ref{finitedual}).
The dimension of an irreducible scheme is the relative dimension over $S$;
in particular, an irreducible scheme of dimension $d$ over
$\eta$ has dimension $d+1$ over $S$, and the complex $\DD_X(n)$
restricted to the generic fiber would be the complex $\DD_{X_\eta}(n-1)[2]$
if it were viewed relative to the generic point.

\begin{theorem}\label{locgeneral}
Let $S$ be a strictly henselian discrete valuation ring of mixed characteristic
with algebraically closed residue field. If $X$ is essentially of finite
type over $S$ and $n\leq 0$, then
$$\DD_X(n)(X)\cong R\Gamma(X_\et,\DD_X(n)).$$
\end{theorem}

\proof
As in the proof of Theorem \ref{locok},
it suffices to show that $\DD(n)(\Spec F)\cong R\Gamma(\Spec F_\et,\DD(n))$
for an extension of finite transcendence degree $d$ of the
residue field of the closed point or the generic point of $S$,
because $\DD_X(n)$ has the localization property for schemes of finite
type over a discrete valuation ring by Levine \cite{levinemoving}.
The case that $F$ lies over the closed point was treated in the proof of
Theorem \ref{locok}, hence we can assume that $F$ lies over $\eta$.
We have to show that the map $H_i(F,\DD(n)))\to H_i(F_\et,\DD(n))$ is an
isomorphism for all $i$ (it is important to remember that we use
dimension relative to $S$ here, so that $F$ is a limit of schemes
of dimension $d+1$). Rationally, Zariski and \'etale
hypercohomology of the motivic complex agree.
With mod $m$-coefficients, we get the isomorphism
from the Beilinson-Lichtenbaum conjecture for $i\geq d+1+n$.
On the other hand, by Proposition \ref{identifyme},
$H_i(F_\et, \DD/m(n))\cong H^{2d+2-i}(F_\et,\mu_m^{\otimes (d+1-n)})$,
and this vanishes for $i<d+1$ because $F$ has cohomological dimension $d+1$.
\proofend

\begin{corollary}\label{setup}
Let $S$ be a Dedekind ring of characteristic $0$ with perfect residue
fields, and let $n\leq 0$.

a) If $i:Z\to X$ is a closed embedding of schemes over $S$ with open
complement $U$, then we have a quasi-isomorphism
$\DD_Z(n)\cong Ri^!\DD_X(n)$, hence a distinguished triangle
$$\cdots \to  R\Gamma(Z_\et,\DD_Z(n))\to R\Gamma(X_\et,\DD_X(n))
\to R\Gamma(U_\et,\DD_U(n))\to \cdots. $$

b) If $f:X\to Y$ is a proper map between schemes over $S$,
we have a push-forward map $f_*:Rf_*\DD_X(n)\to \DD_Y(n)$.
\end{corollary}

\proof
a) Since the statement is local for the \'etale topology,
we can assume that $S$ is strictly henselian local.
Let $s=\Spec k$ be the closed point, $\eta$ be the generic point,
and $X_s$ and $X_\eta$ be the closed and generic fiber, respectively.
We first treat the case $Z=X_s$, $U=X_\eta$.
Consider the following map of distinguished triangles, where the global
section functor is with respect to the \'etale topology,
$$\begin{CD}
\DD(n)(X_s)@>>> \DD(n)(X)@>>> \DD(n)(X_\eta)\\
@VVV @VVV @VVV \\
R\Gamma(X_s,\DD_{X_s}(n))@>>> R\Gamma(X,\DD_X(n))@>>>
R\Gamma(X_\eta,\DD_{X_\eta}(n)).
\end{CD}$$
The vertical maps are quasi-isomorphisms by Theorem \ref{locgeneral},
and the upper triangle is exact by \cite{levinemoving},
hence the lower triangle is exact.
For arbitrary $Z$, we consider the diagram
$$\begin{CD}
R\Gamma(Z_s,\DD(n))@>>> R\Gamma(Z,\DD(n))@>>>
R\Gamma(Z_\eta,\DD(n))\\
@VVV @VVV @VVV \\
R\Gamma(X_s,\DD(n))@>>> R\Gamma(X,\DD(n))@>>>
R\Gamma(X_\eta,\DD(n))\\
@VVV @VVV @VVV \\
R\Gamma(U_s,\DD(n))@>>> R\Gamma(U,\DD(n))@>>>
R\Gamma(U_\eta,\DD(n)).
\end{CD}$$
The horizontal triangles are distinguished by the above, and the outer
vertical triangles are distinguished by Corollary \ref{purity}.
Hence the middle vertical triangle is distinguished as well.
Part b) is proved exactly like Corollary \ref{pushmap},
using Theorem \ref{locgeneral}.
\proofend

\begin{theorem}\label{qumain}
Let $f:X\to S$ be a scheme over a Dedekind ring
of mixed characteristic with perfect residue fields.
Then for every torsion sheaf $\g$ on $X$, and $n\leq 0$,
there is a quasi-isomorphism
$$ R\Hom_X(\g,\DD_X(n))\cong R\Hom_S(Rf_!\g,\DD_S(n)).$$
\end{theorem}

\proof
Replacing $X$ by a compactification $j:X\to T$ and
$\g$ by $j_!\g$, we can assume that $X$ is proper.
Writing $\g$ as a colimit of constructible sheaves, we
can assume that $\g$ is constructible, see the proof
of Corollary \ref{adji}.
The quasi-isomorphism is induced by the map
$Rf_*\DD_X(n)\to \DD_S(n)$ of Corollary \ref{setup}b).
We can assume that $\g$ is $m$-torsion for some integer
$m$, and it suffices to show that the adjoint map
$\DD_X/m(n)\to Rf^!\DD_S/m(n)$ is an isomorphism.
We can check this at stalks and assume that the base
is a henselian discrete valuation ring; the case that the
base is an algebraically closed field is Corollary \ref{adji}.
Consider the following commutative diagram, coming from
Corollary \ref{setup}, and the quasi-isomorphisms
$Rj_*Rf^!=Rf^!Rj_*$ and $i_*Rf^!=Rf^!i_*$:
$$\begin{CD}
i_*\DD_{X_s}/m(n)@>>> \DD_X/m(n)@>>> Rj_*\DD_{X_\eta}/m(n)\\
@VVV @VVV @VVV \\
i_*Rf^!\DD_s/m(n)@>>>Rf^!\DD_S/m(n)@>>> Rj_*Rf^!\DD_\eta/m(n).
\end{CD}$$
The outer maps are quasi-isomorphisms
by Corollary \ref{adji}a), hence so is the middle map.
\proofend

\begin{lemma} \label{strr}
Under the conditions of the Theorem, we have
$\DD_S\cong {\Bbb G}_m[1]$ on $S$.
\end{lemma}

\proof
By \cite[Lemma 11.2]{marcmc}, $\DD_S$ is acyclic (as a complex
of \'etale sheaves) except in degree $-1$.
The quasi-isomorphism is induced by the map
${\Bbb G}_m\to {\mathcal H}^{-1}(\DD_S)$,
sending a unit $u$ on $U$ to the subscheme $(\frac{1}{1-u},\frac{-u}{1-u})$.
\proofend

\subsection{Local duality}
Let $f:X\to S$ be a scheme over a discrete valuation ring and $i:s\to X$
be the closed point. For a torsion
sheaf $\g$, we define cohomology with compact support in the closed
fiber $R\Gamma_{X_s,c}(X_\et,\g)$ to be $R\Gamma(S_\et,i_*Ri^!Rf_!\g)$.

\begin{theorem}\label{localduality}
Let $X\to S$ be a scheme over a henselian discrete
valuation ring of characteristic $0$ with finite residue field. Then
for every torsion sheaf $\g$ on $X$, there is a quasi-isomorphism
$$R\Hom_X(\g,\DD_X)\cong R\Hom_\Ab(R\Gamma_{X_s,c}(X_\et,\g),\Z)[-1].$$
In particular, for constructible $\g$, there are perfect pairings
of finite groups
$$ \Ext^{2-i}_X(\g,\DD_X)\times H^i_{X_s,c}(X_\et,\g)\to \Q/\Z.$$
\end{theorem}

\proof
From Theorem \ref{qumain} we get
$$ R\Hom_X(\g,\DD_X)\cong R\Hom_S(Rf_!\g,\DD_S),$$
hence the result follows from

\begin{lemma} For every complex of constructible sheaves $\f^\cdot$ on $S$, we
have a quasi-isomorphism
$$ R\Hom_S(\f^\cdot,\DD_S)\cong R\Hom_\Ab(R\Gamma_s(S_\et,\f^\cdot),\Z)[-1].$$
\end{lemma}

\proof
This follows using Lemma \ref{strr} from the local duality
quasi-isomorphism \cite[II Thm. 1.8]{adt}
$$ R\Hom_S(\f,\DD_S)[-1]\cong R\Hom_S(\f,{\Bbb G}_m)
\cong R\Hom_\Ab(R\Gamma_s(S_\et,\f),\Q/\Z)[-3].$$
\proofend

Kato-homology $H_i^K(X,\Z/m)$ over a henselian discrete valuation ring
is defined as in \eqref{katocomplex}.
If $X$ is proper and regular over $S$,
with generic fiber $X_\eta$ and closed fiber $X_s$, then
Kato conjectures \cite[Conj. 5.1]{kato} that the Kato-homology
$H_i^K(X_s,\Z/m)$ of the closed and $H_i^K(X_\eta,\Z/m)$
of the generic fiber agree, or
equivalently that the Kato-homology of $X$ vanishes for all $i$.
The same proof as for \eqref{S} gives

\begin{corollary}
For every scheme over a henselian discrete valuation ring of characteristic
$0$ with finite residue fields, there is a long exact sequence
\begin{equation}\label{U}
\cdots \to CH_0(X,i,\Z/m)\to H^{i+1}_{X_s,c}(X_\et,\Z/m)^*
\to H_{i+1}^K(X,\Z/m)\to \cdots .
\end{equation}
\end{corollary}

Note that the exacts sequences \eqref{S}, \eqref{T}, \eqref{U}
fit together into a double-complex
$$\begin{CD}
CH_{-1}(X_\eta,i,\Z/m)@>>> CH_0(X_s,i-1,\Z/m)@>>> CH_0(X,i-1,\Z/m)\\
@VVV @VVV @VdVV \\
H^i_c((X_\eta)_\et,\Z/m)^*@>>>H^i_c((X_s)_\et,\Z/m)^*@>>>
H^i_{X_s,c}(X_\et,\Z/m)^*\\
@VVV @VVV @VVV\\
H_i^K(X_\eta,\Z/m) @>>> H_i^K(X_s,\Z/m) @>>> H_i^K(X,\Z/m)
\end{CD}$$

\subsection{Number rings}
Let $B$ be the spectrum of a number ring. For a torsion sheaf $\f$ on $B$,
the cohomology with compact supports
$R\Gamma_c(B_\et,\f)$ is defined, for example, in \cite[II \S 2]{adt}
and differs from $R\Gamma(B_\et,\f)$ only at the
prime $2$ and only for those $B$ having a real embedding.
For a torsion sheaf $\g$ on $X$, we define cohomology with compact support
$H^i_c(X_\et,\g)$ to be the cohomology of the complex
$R\Gamma_c(B_\et,Rf_!\g)$.
The following generalizes and unifies \cite[II Thms. 6.2, 7.16]{adt}.

\begin{theorem}
For every scheme $f:X\to B$ and torsion sheaf $\g$ on $X$,
we have a quasi-isomorphism
$$R\Hom_X(\g,\DD_X)\cong R\Hom_\Ab(R\Gamma_c(X_\et,\g),\Z)[-1],$$
which induce perfect pairings of finite groups for constructible $\g$,
$$ \Ext^{2-i}_X(\g,\DD_X)\times H^i_c(X_\et,\g)\to \Q/\Z.$$
\end{theorem}

\proof
For a complex of constructible sheaves $\f^\cdot$ on $B$, we have
by Artin-Verdier duality \cite{mazur}\cite[II Thm.3.1b)]{adt}
a quasi-isomorphism
$$ R\Hom_B(\f^\cdot, \G_m)\cong R\Hom_\Ab(R\Gamma_c(B,\f^\cdot),\Q/\Z)[-3].$$
If we apply this to the quasi-isomorphism of Theorem \ref{qumain}, we get
\begin{multline*}
R\Hom_X(\g,\DD_X)\cong R\Hom_B(Rf_!\g,\G_m)[1]\\
\cong R\Hom_\Ab(R\Gamma_c(X_\et,\g),\Q/\Z)[-2].
\end{multline*}
\proofend


Over the spectrum of a number ring, higher Chow groups
$CH_i(X,n)$ are defined as the Zariski-hypercohomology of
the complex $z_n(-,*)$, and Kato-homology $H_i^K(X,\Z/m)$ is
defined \cite[Conj. 05]{kato}
as the homology of the cone of the canonical map of
the complex \eqref{katocomplex} to the direct sum of the complexes
\eqref{katolocal} for $X\times_BF_\nu$, where $F_\nu$ runs through
the real places of $B$. Kato conjectures that $H_i^K(X,\Z/m)=0$
for $i>0$ and $X$ regular, flat and proper over $B$.

\begin{corollary}
For every scheme over a number ring, there is a long exact sequence
$$ \cdots \to CH_0(X,i,\Z/m)
\to H^{i+1}_c(X_\et, \Z/m)^* \to H_{i+1}^K(X,\Z/m)\to \cdots.$$
\end{corollary}

For completeness we give the following analog of Corollary \ref{adji}
and Theorem \ref{exchange1}

\begin{proposition}\label{adjiglobal}
Let $f:X\to Y$ be a map of schemes over the spectrum of a number ring $S$,
and let $n\leq 0$.

a) For every locally constant constructible sheaf $\f$ on $Y$, the map of
Corollary \ref{setup}b) induces a quasi-isomorphism
$$ \DD_X(n)\otimes f^*\f\cong Rf^!(\DD_Y(n)\otimes\f) .$$

b) For every torsion sheaf $\g$ on $X$ and every finitely generated,
locally constant sheaf $\f$ on $Y$, we have a functorial quasi-isomorphism
$$ \Hom_Y(Rf_!\g,\DD_Y(n)\otimes\f)\cong \Hom_X(\g,\DD_X(n)\otimes f^*\f).$$

c) (Exchange formulas)
If $\g$ and $\f$ are constructible, then
\begin{align*}
\du(\g\otimes \f)&\cong R\HOM(\g,\du(\f))\\
Rf_*\du(\f) &\cong \du (Rf_!\f) \\
Rf^!\du(\f) &\cong \du (f^*\f).
\end{align*}
\end{proposition}

\proof
a) Since the statement is local for the \'etale topology, we can assume
that $S$ is strictly henselian
and that $\f$ is a constant sheaf of the form $\Z/m$.
Then the proof of Theorem \ref{qumain} works in this situation.
b) follows from a) as in Corollary \ref{adji}, and
c) is proved as in Theorem \ref{exchange1}.
\proofend

\end{document}